\documentclass[preprint,12pt]{elsarticle}
\usepackage{stmaryrd}
\usepackage{mathrsfs}
\usepackage{amsmath}
\usepackage{amsfonts}

\usepackage{amssymb}
\usepackage{amsthm}
\usepackage{amsmath}

\journal{\quad }

\begin{document}

\begin{frontmatter}



\title{On the Rayleigh--Taylor instability for the incompressible\\
 viscous magnetohydrodynamic equations}

\author[FJ]{Fei Jiang\corref{cor1}}
\ead{jiangfei0591@163.com}
\author[FJ]{Song Jiang}
 \ead{jiang@iapcm.ac.cn} \cortext[cor1]{Corresponding author: Tel +86
15001201710.}
\author[YJ]{Yanjin Wang}
\ead{wangyanjin\mbox{\underline{\ }}2008@163.com}
\address[FJ]{Institute of Applied Physics and Computational Mathematics, Beijing, 100088, China.}
\address[YJ]{School of Mathematical Sciences, Xiamen University, Fujian 361005, China.}

\begin{abstract} We study the Rayleigh-Taylor instability
problem for two incompressible, immiscible, viscous
magnetohydrodynamic (MHD) flows with zero resistivity and surface
tension (or without surface tension), evolving with a free interface
under presence of a uniform gravitational field. First, we
reformulate the MHD free boundary problem in a infinite slab as a
Navier-Stokes system in Lagrangian coordinates with a force term
induced by the fluid flow map. Then, we analyze the linearized
problem around the steady state which describes a denser immiscible
fluid lying above a light one with a free interface separating the
two fluids, and both fluids being in (unstable) equilibrium. By
studying a family of modified variational problems, we construct
smooth (when restricted to each fluid domain) solutions to the
linearized problem that grow exponentially fast in time in Sobolev
spaces, thus leading to an global instability result for the
linearized problem. Finally, using these pathological solutions, we
prove the global instability for the corresponding nonlinear problem
in an appropriate sense. Moreover, we evaluate that the so-called
critical number indeed is equal to $\sqrt{g[\varrho]/2}$, and
analyze the effect of viscosity and surface tension on the
instability.

\end{abstract}

\end{frontmatter}


\newtheorem{thm}{Theorem}[section]
\newtheorem{lem}{Lemma}[section]
\newtheorem{pro}{Proposition}[section]
\newtheorem{cor}{Corollary}[section]
\newproof{pf}{Proof}
\newdefinition{rem}{Remark}[section]
\newtheorem{definition}{Definition}[section]
\section{Introduction}
\label{Intro} \numberwithin{equation}{section} Considering two
completely plane-parallel layers of immiscible fluids, the heavier
on top of the lighter one and both subject to the earth's gravity.
In this case, the equilibrium state is unstable to sustain small
perturbations or disturbances, and this unstable disturbance will
grow and lead to a release of potential energy, as the heavier fluid
moves down under the (effective) gravitational force, and the
lighter one is displaced upwards. This phenomenon was first studied
by Rayleigh \cite{RLAP,RLIS} and then Taylor \cite{TGTP}, and is
called therefore the Rayleigh-Taylor (R-T) instability. In the last
decades, this phenomenon has been extensively investigated from both
physical and numerical aspects and there are a lot of related
results in the literature. In particular, many results concerning
the linearized problems have been summarized in monographs, see, for
instance, \cite{CSHHS,WJH}. To our best knowledge, however, there
are only few mathematical analysis results on nonlinear problems in
the literature, due to the fact that in general, passage from a
linearized instability to a dynamical nonlinear instability for a
conservative nonlinear partial differential system is rather
difficult.

The magnetohydrodynamic (MHD) analogue of the R-T
instability arises when the fluids are electrically conducting and a
magnetic field is present, and the growth of the instability will be
influenced by the magnetic field due to the generated
electromagnetic induction and the Lorentz force.
 This has been
analyzed from the physical point of view in many monographs, see,
for example, \cite{CSHHS,WJH}. Because of additional difficulties
induced by presence of the magnetic field, many results concerning
the R-T instability of superposed flows could not be directly
generalized to the case of MHD flows.

 In this paper, we study the R-T problem
for two incompressible, immiscible, viscous magnetohydrodynamic
(MHD) flows, with zero resistivity and surface tension (or without
surface tension), evolving with a free interface in the presence of
a uniform gravitational field. We will prove that in Lagrangian
coordinates, the corresponding MHD linearized system is globally
unstable in some sense  as time increases, and moreover, the
original nonlinear problem with or without surface tension is
globally unstable in some appropriate sense.  Next, we formulate our
problem in details.

\subsection{Formulation in Eulerian coordinates}
We consider the two-fluids free boundary problem for the equations of
magnetohydrodynamics (MHD) within the infinite slab
$\Omega=\mathbb{R}^{2}\times (-1,1)\subset \mathbb{R}^3$ and for
time $t\geq 0$. The fluids are separated by a moving free interface
$\Sigma (t)$ that extends to infinity in every horizontal direction.
The interface divides $\Omega$ into two time-dependent, disjoint,
open subsets $\Omega_{\pm}(t)$ so that $\Omega=\Omega_+(t)\cup
\Omega_-(t)\cup \Sigma (t)$ and $\Sigma (t)=\bar{\Omega}_+(t)\cap
\bar{\Omega}_-(t)$. The motions of the fluids are driven by the
constant gravitational field along $e_3$ (the $x_3$-direction),
$G=(0,0,-g)$ with $g>0$ and the Lorentz force induced by the
magnetic fields. The motion of the fluids is described by their velocity,
pressure and magnetic field functions, which are given for each
$t\geq 0$ by, respectively,
\begin{equation*}\label{0101}
(u_\pm,\bar{p}_\pm, h_\pm)(t,\cdot):\Omega_\pm(t)\rightarrow
(\mathbb{R}^3,\mathbb{R}^+,\mathbb{R}^3).
\end{equation*}
We assume that at a given time $t\geq 0$, these functions have
well-defined trace onto $\Sigma(t)$.

The fluids under consideration are incompressible, viscous and of
zero resistivity. Hence for $t>0$ and $x=(x_1,x_2,x_3)\in \Omega_\pm(t)$,
the fluids satisfy the following magnetohydrodynamic equations:
\begin{equation}\label{0102}\left\{
                  \begin{array}{ll}
\partial_t(\varrho_\pm u_\pm)+\mathrm{div}(\varrho_\pm
u_\pm\otimes u_\pm)+\mathrm{div}S_\pm=-g\varrho_\pm e_3,  \\[1mm]
                    \mathrm{div}u_\pm=0,  \\[1mm]
                    \partial_t h_\pm+\mathrm{div}(u_\pm\otimes
h_\pm)-\mathrm{div}(h_\pm\otimes u_\pm)=0,\\[1mm]
\mathrm{div}h_\pm=0,
                  \end{array}
                \right.
\end{equation}
where we have defined the stress tensor consisting of both fluid and
magnetic parts by
\begin{equation*}\label{0103}
S_\pm=-\mu_\pm(\nabla u_\pm+\nabla u_\pm^T )+\bar{p}_\pm I+\frac{|h_\pm|^2}{2}I-h_\pm\otimes h_\pm.
\end{equation*}
Hereafter the superscribe $T$ means the transposition and $I$ is the
$3\times 3$ identity matrix. The positive constants $\varrho_\pm$
denote the densities of the respective fluids.

For two viscous fluids meeting at a free boundary, the standard
assumptions are  that the velocity is continuous across the
interface and that the jump in the normal stress is proportional to
the mean curvature of the surface multiplied by the normal vector to
the surface (cf. \cite{CSHHS,WJLE}). This requires us to enforce the
jump conditions
\begin{eqnarray}
&&\label{0104}[u]|_{\sum (t)}=0,\\[0.5em]
&&\label{0105}[S \nu]|_{\sum (t)}=\kappa H\nu,
\end{eqnarray}
where we have written the normal vector to $\Sigma(t)$ as $\nu$,
$f|_{\Sigma(t)}$ for the trace of a quantity $f$ on $\Sigma(t)$, and
denoted the interfacial jump by
\begin{equation*}\label{0106}[f]|_{\Sigma(t)}:=f_+|_{\Sigma(t)}-f_-|_{\Sigma(t)}.
\end{equation*}
Here we take $H$ to be twice the mean curvature of the surface
$\Sigma(t)$ and the surface tension to be a constant $\kappa\geq 0$.
We will also enforce the condition that the fluid velocity vanishes
at the fixed boundaries; we implement this via the boundary
conditions
\begin{equation*}\label{0107}
u_+(t,x',-1)=0,\quad u_-(t,x',1)=0\;\;\mbox{ for all }t\geq 0,\
x'=(x_1,x_2)\in \mathbb{R}^{2}.
\end{equation*}

Since the fluids are of zero resistivity, the magnetic equations
(\ref{0102})$_3$ are a free transport system along the flow, and
hence the Dirichlet boundary condition on the velocity at the fixed
boundary prevents the necessity of prescribing boundary condition on
the magnetic field. On the other hand, due to the incompressibility
(\ref{0102})$_4$ and also from the physical point of view, we assume
that the normal component of the magnetic field is continuous across the
free interface (cf. \cite{CSHHS,SJA})
\begin{equation}\label{0108}[h\cdot \nu]|_{\sum(t)}=0.
\end{equation}
In fact, we will show in the next subsection that the
incompressibility of $h_\pm$ and the jump condition (\ref{0108}) are satisfied
if they hold initially. Therefore, the conditions
(\ref{0102})$_4$, (\ref{0108}) are transformed to the compatibility
conditions assumed on the initial magnetic field.

 The motion of the free interface is coupled to the evolution
equations for the fluids (\ref{0102}) by requiring that the surface
be advected with the fluids. More precisely, if $V(t,x)\in
\mathbb{R}^3$ denotes the normal velocity of the boundary at $x\in
\Sigma(t)$, then
\begin{equation*}\label{0109} V(t,x)=\big( u(t,x)\cdot \nu(t,x)\big)\nu(t,x) ,
\end{equation*}
where $u(t,x)$ is the common trace of $u_\pm(t,x)$ onto $\Sigma(t)$
and these traces agree because of the jump condition (\ref{0104}),
which also implies that there is no possibility of the fluids
slipping past each other along $\Sigma(t)$.

To complete the statement of the problem, we must specify initial
conditions. We give the initial interface $\Sigma(0)=\Sigma_0$
(e.g., $\Sigma_0=\{x\in {\mathbb R}^3\,|\; x_3=0\}$),
which yields the open sets $\Omega_\pm(0)$ on which we
specify the initial data for the velocity and magnetic field
\begin{equation*}(u_\pm,h_\pm)(0,\cdot):\Omega_\pm(0)\rightarrow
(\mathbb{R}^3,\mathbb{R}^3).
\end{equation*}
Thus the initial datum of the pressure $\bar{p}_0$ can be defined by
$\varrho_\pm$, $\Sigma_0$,  $u_\pm(0,\cdot)$ and $h_\pm(0,\cdot)$.

 To simply the equations we introduce the indicator
function $\chi$ and denote
\begin{equation*}
            \begin{array}{ll}
           \varrho=\varrho_+\chi_{\Omega_+}+\varrho_-\chi_{\Omega_-},      &
   u=u_+\chi_{\Omega_+}+u_-\chi_{\Omega_-},\\[0.5em]
h=h_+\chi_{\Omega_+}+h_-\chi_{\Omega_-}, &
\bar{p}=\bar{p}_+\chi_{\Omega_+}+\bar{p}_-\chi_{\Omega_-},
                  \end{array}
\end{equation*}
and also define the modified pressure by
\begin{equation*}\label{0110}
p=\bar{p}+\frac{|h|^2}{2}+g\varrho x_3.
\end{equation*}
Hence, the equations (\ref{0102}) can be rewritten as
\begin{equation*}\label{0111}
\left\{
            \begin{array}{l}
                \varrho\partial_t u+\varrho u\cdot\nabla u+\nabla p=h\cdot \nabla h+\mu \Delta u,
\\
\partial_th+u\cdot \nabla h-h\cdot\nabla u=0,\\
\mathrm{div}u=\mathrm{div}h=0
                  \end{array}    \right.
\end{equation*}
in $\Omega\backslash\Sigma (t)$ for each $t > 0$, and the jump condition (\ref{0105}) becomes,
setting $[\varrho]=\varrho_+-\varrho_-$,
\begin{equation*}\label{0112}
[(pI-\mu(\nabla u+\nabla
u^T))\nu]\big|_{\sum(t)}=g[\varrho]x_3\nu+h\cdot\nu[h]\big|_{\sum(t)}+\kappa H\nu.
\end{equation*}

\subsection{Formulation in Lagrangian coordinates}
The movement of the free interface $\Sigma(t)$ and the subsequent change of the
domains $\Omega_\pm(t)$ in Eulerian coordinates will result in severe
mathematical difficulties. To circumvent such difficulties, 
we switch our analysis to Lagrangian
coordinates, so that the interface and the domains stay fixed in time. To
this end, we define the fixed Lagrangian domains
$\Omega_+=\mathbb{R}^{2}\times (0,1)$ and
$\Omega_-=\mathbb{R}^{2}\times (-1,0)$, and assume that there exist
invertible mappings
\begin{equation*}\label{0113}
\eta_\pm^0:\Omega_\pm\rightarrow \Omega_\pm(0),
\end{equation*}
such that $\Sigma_0=\eta_+^0(\{x_3=0\})$,
$\{x_3=1\}=\eta_+^0(\{x_3=1\})$ and
$\{x_3=-1\}=\eta_-^0(\{x_3=-1\})$. The first condition means that
$\Sigma_0$ is parameterized by the mapping $\eta_+^0$
restricted to $\{x_3=0\}$, while the latter two conditions mean that
$\eta_\pm^0$ map the fixed upper and lower boundaries into
themselves. Define the flow maps $\eta_\pm$ as the solution to
\begin{equation}\label{0114}
\left\{
            \begin{array}{l}
\partial_t \eta_\pm(t,x)=u_\pm(t,\eta_\pm(t,x))
\\
\eta_\pm(0,x)=\eta_\pm^0(x).
                  \end{array}    \right.
\end{equation}

We denote the Eulerian coordinates by $(t,y)$ with $y=\eta(t,x)$,
whereas the fixed $(t,x)\in \mathbb{R}^+\times \Omega$ stand for the
Lagrangian coordinates. In order to switch back and forth from
Lagrangian to Eulerian coordinates, we assume that
$\eta_\pm(t,\cdot)$ are invertible and
$\Omega_{\pm}(t)=\eta_{\pm}(t,\Omega_{\pm})$, and since $u_\pm$ and
$\eta_\pm^0$ are all continuous across $\{x_3=0\}$, we have
$\Sigma(t)=\eta_\pm(t,\{x_3=0\})$. In other words, the Eulerian
domains of upper and lower fluids are the image of $\Omega_\pm$
under the mappings $\eta_\pm$, and the free interface is the image
of $\{x_3=0\}$ under the mappings $\eta_\pm(t,\cdot)$.

Setting $\eta=\chi_+\eta_++\chi_-\eta_-$, we define the Lagrangian
unknowns by
\begin{equation*}\label{0115}
(v,q,b)(t,x)=(u,p,h)(t,\eta(t,x)),\ (t,x)\in \mathbb{R}^+\times
\Omega.
\end{equation*}
Defining the matrix $A:=(A_{ij})_{3\times 3}$ via
$A^T=(D\eta)^{-1}:=(\partial_j \eta_i)^{-1}_{3\times 3}$, and the
identity matrix $I=(I_{ij})_{3\times 3}$,  thus in Lagrangian
coordinates the evolution equations for $\eta$, $v$, $q$, $b$ read as
(writing $\partial_j=\partial/\partial_{x_j}$)
\begin{equation}\label{0116} \left\{
                              \begin{array}{ll}
\partial_t \eta_i=v_i \\[1mm]
\varrho\partial_t v_i+A_{jk}\partial_k T_{ij} =b_jA_{jk}\partial_k b_i, \\[1mm]
A_{jk}\partial_{k}v_j=0, \\[1mm]
\partial_t b_i=b_jA_{jk}\partial_kv_i, \\[1mm]
A_{jk}\partial_k b_j=0,
\end{array}
                            \right.
\end{equation}
where the stress tensor of fluid part in Lagrangian coordinates,
$T(v,q)$, is given by
\begin{equation*}\label{0117}T_{ij}=q I_{ij}-\mu(A_{jk}\partial_k v_i+A_{ik}\partial_k
v_j).
\end{equation*}
Here we have written $I_{ij}$ for the $ij$-component of the identity
matrix $I$ and used the Einstein convention of summing over repeated indices.

To write the jump conditions, for a quantity $f=f_{\pm}$, we define
the interfacial jump by
\begin{equation*}\llbracket  f
\rrbracket:=f_+|_{x_3=0}-f_-|_{x_3=0}.
\end{equation*}
Then the jump conditions in Lagrangian coordinates are
\begin{equation}\label{0118}
\llbracket  v \rrbracket =0,\ \llbracket  b_j n_j \rrbracket=0,\
\llbracket  T_{ij}n_j \rrbracket=g[\varrho]\eta_3 n_i+\llbracket b_i
\rrbracket b_j n_j+\kappa H n_i,\end{equation} where we have written
$n:=(n_1, n_2, n_3)=\nu(\eta)$, i.e.
\begin{equation}\label{0119}n=\frac{\partial_1\eta\times \partial_2\eta}{|\partial_1\eta\times
\partial_2\eta|}=\frac{Ae_3}{|Ae_3|}\bigg|_{\{x_3=0\}}
\end{equation}
for the normal vector to the surface $\Sigma(t)=\eta(t,\{x_3=0\})$
and $H$ for twice the mean curvature of $\Sigma(t)$. Since
$\Sigma(t)$ is parameterized by $\eta$, we may employ the standard
formula for the mean curvature of a parameterized surface to write
\begin{equation}\label{0120}H=\left(\frac{|\partial_1\eta|^2
\partial_2^2\eta-2(\partial_1\eta\cdot\partial_2\eta)\partial_1\partial_2\eta+
|\partial_2\eta|^2\partial_1^2
\eta}{|\partial_1\eta|^2|\partial_2\eta|^2-|\partial_1\eta\cdot\partial_2\eta|^2}\right)\cdot
n.
\end{equation}
Finally, we require the no-slip boundary conditions
\begin{equation}\label{0121}
v_-(t,x',-1)=0,\quad v_+(t,x',1)=0\;\;\mbox{ for all }t\geq 0,\ x'\in\mathbb{R}^2.
\end{equation}

\subsection{Reformulation}
In this subsection we reformulate the free boundary problem
(\ref{0116}), (\ref{0118}) and (\ref{0121}). Our goal is to
eliminate $b$ by expressing it in terms of $\eta$, and this can be
achieved in the same manner as  in \cite{WYC}. For the reader's
convenience, we give the derivation here.

Applying $A_{il}$ to (\ref{0116}$)_4$, we  have
\begin{equation*}\label{0122}A_{il}\partial_t b_i=b_jA_{jk}\partial_k v_i
A_{il}=b_jA_{jk}\partial_t(\partial_k
\eta_i)A_{il}=-b_jA_{jk}\partial_k \eta_i \partial_t
A_{il}=-b_i\partial_tA_{il},
\end{equation*}
which implies that $\partial_t(A_{jl}b_j)=0$. Hence,
\begin{eqnarray}  \label{0123}   A_{jl}b_j=A^0_{jl}b_j^0, \\[1mm]
\label{0124}   b_i=\partial_l \eta_i A^0_{jl}b_j^0.
\end{eqnarray}
Hereafter, the superscript $0$ means the initial value.

With help of (\ref{0124}), and  the geometric identities
$J=J^0$, $\partial_k(JA_{ik})=0$ where $J=|D\eta|$,
we can apply $A_{ik}\partial_k$ to (\ref{0124}) to deduce that the divergence of $b$
(i.e. (\ref{0116}$)_5$) satisfies
\begin{equation}\label{0126}A_{ik}\partial_k b_i=\frac{J}{J^0}A_{ik}\partial_k(\partial_l \eta_i A_{jl}^0
b_j^0)=\frac{1}{J^0}\partial_k(JA_{ik}\partial_l \eta_i A_{jl}^0
b^0_j)=\frac{1}{J^0}\partial_k(J^0A_{jk}^0 b_j^0)=A^0_{jk}\partial_k
b_j^0.
\end{equation}

Next, we evaluate the jump $\llbracket  b_j n_j \rrbracket $. It is easy,
recalling (\ref{0119}), to verify that $Ae_3$ is continuous across
the free interface. Hence, one has
\begin{equation}\label{0127}
\llbracket b_j n_j \rrbracket=\llbracket \partial_l \eta_j A_{kl}^0b_k^0  A_{j3}
\rrbracket \frac{1}{|Ae_3|}=\llbracket A_{k3}^0 b_k^0
\rrbracket\frac{1}{|Ae_3|}=\llbracket  b_j^0 n_j^0
\rrbracket\frac{|A^0e_3|}{|Ae_3|}.
\end{equation}
Thus, if one assumes the compatibility conditions on the initial data
\begin{equation}\label{0128} A_{jk}^0\partial_k b_j^0=0,\;\;\;\llbracket  b_j^0 n_j^0
\rrbracket=0,
\end{equation}
then from (\ref{0126}) and (\ref{0127}) one gets
\begin{equation}\label{0129} A_{jk}\partial_k b_j=0,\;\;\;\llbracket  b_j n_j \rrbracket=0.
\end{equation}

Moreover, for simplifying notation, we assume that
\begin{equation}\label{0130} A_{ml}^0b_m^0=\bar{B_l}\;\;\mbox{ with }\bar{B}\mbox{ being a
constant vector}.
\end{equation}
We remark that the class of the pairs of the data $\eta^0$, $b^0$
that satisfy the constraints (\ref{0128}) and  (\ref{0130}) is quite
large. For example, we chose $\eta^0=$Id and $b^0=$constant vector,
then $\eta^0$ and $b^0$ satisfy (\ref{0128}) and  (\ref{0130}), consequently,
the pair ($\eta$, $b$) which is transported by
the flow will satisfy (\ref{0129}).

Now, in view of (\ref{0123}), (\ref{0124}), (\ref{0130}), we
represent the Lorentz force term by
\begin{equation*}\label{0131}b_j A_{jk}\partial_k b_i=\partial_l \eta_j A_{ml}^0
b^0_m A_{jk}\partial_k(\partial_r \eta_i A_{sr}^0b_s^0)=A_{mk}^0
b_m^0\partial_k(\partial_r \eta_i A_{sr}^0 b_s^0)=\bar{B}_l\bar{B}_m
\partial_{lm}^2\eta_i.\end{equation*}
Hence the equations (\ref{0116}) become a Navier-Stokes system
with the force term induced by the flow map $\eta$:
\begin{equation}\label{0132} \left\{
                              \begin{array}{ll}
\partial_t \eta_i=v_i\\[1mm]
\varrho\partial_t v_i+A_{jk}\partial_k
T_{ij}=\bar{B}_l\bar{B}_m\partial_{lm}^2 \eta_i,\\[1mm]
A_{jk}\partial_k v_j=0,
\end{array}
                            \right.
\end{equation}
where the magnetic number $\bar{B}$ can be regarded as a vector
parameter. Accordingly, the jump conditions (\ref{0118}) become
\begin{equation}\label{0133}\begin{aligned}
\llbracket  v \rrbracket =0,\;\;\;\llbracket  T_{ij}n_j
\rrbracket=g[\varrho]\eta_3 n_i+\bar{B}_l\bar{B}_m\llbracket
\partial_l \eta_i \rrbracket \partial_m\eta_j n_j+\kappa H n_i.
\end{aligned}\end{equation}
Note that we implicitly admit that $\bar{B}_m\partial_m \eta_j n_j$
is continuous across $\{x_3=0\}$. In fact, this follows from the assumptions
(\ref{0128}), (\ref{0130}).

Finally, we require the boundary conditions (\ref{0121}).

\subsection{Linearization around the steady state}
The system (\ref{0132}), (\ref{0133}) and (\ref{0121}) admits the
steady solution with $v=0$, $\eta=\mathrm{Id}$, $q=$constant with
the interface given by $\eta(\{x_3=0\})=\{x_3=0\}$, and hence, $n=e_3$
and $A=I$. Now we lnearize the equations (\ref{0132}) around such a
steady-state solution, the resulting linearized equations are
\begin{equation}\label{0134} \left\{
                              \begin{array}{l}
\partial_t\eta=v,\\
\varrho \partial_t v+\nabla q=\mu \Delta v+\bar{B}_l\bar{B}_m\partial^2_{lm}\eta,\\
\mathrm{div} v=0.
\end{array}
                            \right.
\end{equation}
 The corresponding linearized jump conditions are
\begin{equation}\label{0135}
\llbracket  v \rrbracket=0,\;\;\;\llbracket qI-\mu(\nabla v+\nabla
v^T)\rrbracket e_3=(g[\varrho]\eta_3+\kappa \Delta_{x'}\eta_3 )e_3
+\bar{B}_3\bar{B}_l\llbracket
\partial_l \eta \rrbracket,
\end{equation}
while the boundary conditions are
\begin{equation}\label{0136}
v_-(t,x',-1)=0,\quad v_+(t,x',1)=0.
\end{equation}

Now, we briefly review some of the previous results on the nonlinear
 Rayleigh-Taylor instability without the magnetic field.
In 1987, \citet{EDGTC1111} proved the ill-posedness of the equations
of motion for a perfect fluid with free boundary. Then, he adapted
the approach of \cite{EDGTC1111} to obtain the nonlinear
ill-posedness of both R-T and Helmholtz instability problems for
two-dimensional incompressible, immiscible, inviscid fluids without
surface tension \cite{EDGIC1111}. In 2011, for two-compressible
immiscible fluids evolving with a free interface, Guo and Tice made
use of flow maps (cf. (\ref{0114})) to transfer the free boundary
into a fixed boundary and established a variational framework for
the nonlinear instability in \cite{GYTI1}, in which with the help of
the method of Fourier synthesis, they constructed solutions that
grow arbitrarily fast in time in Sobolev spaces, leading to the
ill-posedness of the perturbed problem in Lagrangian coordinates.
This is in some sense a compressible analogue to the ill-posedness
of the R-T instability problem for incompressible fluids given in
\cite{EDGTC1111}. Unfortunately, the approaches in both
\cite{EDGIC1111} and \cite{GYTI1} could not be directly applied to
the viscous flow case, since the viscosity can bring some technical
difficulties in the study of the nonlinear R-T instability. Hence,
Guo and Tice only investigated the stabilized effect of viscosity
and surface tension to the linear R-T instability (see \cite{GYTI2})
for compressible flows, and the corresponding nonlinear instability
still remains open. Recently, Jiang, Jiang and Wang \cite{JFJSWWWO}
showed the nonlinear instability in some sense for two
incompressible immiscible fluids with or without surface tension in
Eulerian coordinates.

Unfortunately,  the results concerning the nonlinear R-T instability
in \cite{EDGIC1111,GYTI1,JFJSWWWO} are still not generalized to the
case of MHD flows, because of additional difficulties induced by
presence of the magnetic field. For two incompressible immiscible
fluids evolving with a free interface (the density is discontinuous
across the interface) under presence of a magnetic field, it was
first shown by Kruskal and Schwarzschild \cite{KMSMSP} that a
horizontal magnetic field has no effect on the development of the
R-T instability for the linearized equations in the whole space in
Eulerian coordinates. Recently, for the case of finite slab domain
without surface tension, Wang \cite{WYC} obtained the critical
magnetic number (denoted by $|B|_c$) for the linear R-T instability.
Namely, he gave an instability criterion to the linearized problem
(\ref{0134}) in terms of the value of the magnetic number, and
pointed out the linear R-T instability in the case $\bar{B}=(0,0,B)$
with $|B|<|B|_c$. In particular, unlike 2D case, he also showed that
the linearized problem is unstable for any initially horizonal
magnetic field $\bar{B}=(B,0,0)$. Form \cite{WYC}, we see
that the stabilized effect of the magnetic field is more
remarkable than that of the surface tension $\kappa$ which only
stabilizes the frequency interval $(0,\sqrt{g[\varrho]/\kappa})$ for
any $\kappa>0$, also see \cite{GYTI2}.  Adopting and modifying the
approaches in both \cite{GYTI1} and \cite{WYC}, Duan, Jiang and
Jiang recently showed the ill-posedness of the linearized system
(\ref{0134}) with $\mu=0$ and $\bar{B}=(B,0,0)$ in the sense of
Hadamard \cite{DRJFJS}, and moreover, the ill-posedness of the
corresponding nonlinear problem in some sense. Finally, we mention
that Hwang \cite{HHVQ} derived the nonlinear instability around
different steady states for both incompressible and compressible
inviscid MHD flows with continuous density, thus extending the
results in \cite{HHJGY} to the case of MHD flows with continuous
density.

 In this paper, we will study the R-T instability for the problem (\ref{0132}), (\ref{0133}) and
(\ref{0121}). We will prove that the corresponding linearized
system (\ref{0134})--(\ref{0136}) with any $\bar{B}=(B,0,0)$ is
globally unstable. Moreover, the original nonlinear problem
(\ref{0132}) with or without surface tension is globally unstable in
some sense. For this purpose, we assume that $\kappa\geq 0$ and that
the upper fluid is heavier than the lower fluid, i.e.,
\begin{equation*}\label{}
\varrho_+>\varrho_-\Leftrightarrow [\varrho]>0.
\end{equation*}  In
addition, we will compute that the so-called critical number in
\cite{WYC} indeed is equal $\sqrt{g[\varrho]/2}$ (see Remark
\ref{rem:03011111} for the details), and analyze the effect
of viscosity and surface tension on the linear instability (see
Remark \ref{rem:0302111111}). The results of the current paper
extend the ones in \cite{HHVQ,WYC}.

 The rest of this paper is organized as follows. In
Section 2 we state our results concerning the linearized system
(\ref{0134}) and nonlinear system (\ref{0132}), i.e. Theorems
\ref{thm:0201} and \ref{thm:0202}, respectively. In Section 3 we
construct the growing solutions to the linearized equations, while
in Section 4 we prove the uniqueness of the linearized problem and
Theorem \ref{thm:0201}. In Section 5, we prove the global
instability of order $k$ to the corresponding nonlinear problem,
i.e., Theorem \ref{thm:0202}.

\section{Main results}

Before stating the main results, we introduce the notation that will
be used throughout the paper. For a function $f\in L^2(\Omega)$, we
define the horizontal Fourier transform via
\begin{equation*}\label{0201}\hat{f}(\xi,x_3)=\int_{\mathbb{R}^2}f(x',x_3)e^{-ix'\cdot
\xi}\mathrm{d}x',
\end{equation*}
where $x', \xi\in \mathbb{R}^2$ and ``$\cdot$" denotes scalar
product. By the Fubini and Parseval theorems, we have
\begin{equation*}\label{0202}
\int_\Omega|f(x)|^2\mathrm{d}x=
\frac{1}{4\pi^2}\int_\Omega\left|\hat{f}(\xi,x_3)\right|^2\mathrm{d}\xi\mathrm{d}x_3.
\end{equation*}

We now define a function space suitable for our analysis  of two
disjoint fluids. For a function $f$ defined on $\Omega$, we write
$f_+$ for the restriction to $\Omega_+=\mathbb{R}^2\times (0,1)$ and
$f_-$ for the restriction to $\Omega_-=\mathbb{R}^2\times (-1,0)$.
For $s\in \mathbb{R}$, define the piecewise Sobolev space of order
$s$ by
\begin{equation}\label{0203}
\bar{\mathrm{H}}^s(\Omega)=\{f~|~f_+\in H^s(\Omega_+),f_-\in
H^s(\Omega_-)\}
\end{equation}
endowed with the norm
$\|f\|_{\bar{\mathrm{H}}^s}^2=\|f\|_{H^s(\Omega_+)}^2+\|f\|^2_{H^s(\Omega_-)}$.
For $k\in \mathbb{N}$ we can take the $\bar{\mathrm{H}}^s$-norm to be
\begin{equation*}\label{0204}\begin{aligned}
\|f\|^2_{\bar{\mathrm{H}}^k}:=&\sum_{j=0}^k\int_{\mathbb{R}^2\times
I_\pm}(1+|\xi|^2)^{k-j}\left|\partial_{x_3}^j\hat{f}_\pm(\xi,x_3)\right|^2\mathrm{d}\xi\mathrm{d}x_3\\
=
&\sum_{j=0}^k\int_{\mathbb{R}^2}(1+|\xi|^2)^{k-j}\left\|\partial_{x_3}^j
\hat{f}_\pm(\xi,\cdot)\right\|^2_{L^2(I_\pm)}\mathrm{d}\xi
\end{aligned}\end{equation*}
 for $I_-=(-1,0)$ and $I_+=(0,1)$. The main difference between the
piecewise Sobolev space $\bar{\mathrm{H}}^s(\Omega)$ and the usual
Sobolev space $H^s(\Omega)$ is that we do not require
functions in the piecewise Sobolev space to have weak derivatives
across the interface $\{x_3=0\}$.

 Now, we are in a position to state our first result, i.e. the
result of global instability for the linearized problem (\ref{0134}).
\begin{thm}\label{thm:0201}
Assume the constant vector $\bar{B}=(B, 0,0)$,  then the linearized
problem (\ref{0134}) with the corresponding jump and boundary
conditions (\ref{0135}), (\ref{0136}) is globally unstable in
$\bar{\mathrm{H}}^k(\Omega)$ for every $k$.  More precisely, there
exists a constant $C_1>0$, and for any $k$, $j\in \mathbb{N}$ with
$j\geq k$ and for any $\alpha>0$, there exist a constant $C_{j,k}$
depending on $j$ and $k$, and a sequence of solutions
$\{(\eta_n,v_n,q_n)\}_{n=1}^\infty$ to (\ref{0134}) satisfying the
corresponding jump and boundary conditions (\ref{0135}),
(\ref{0136}), such that
\begin{equation}\label{0205}
\|\eta_n(0)\|_{\bar{\mathrm{H}}^j}+\|v_n(0)\|_{\bar{\mathrm{H}}^j}+\|q_n(0)\|_{\bar{\mathrm{H}}^j}\leq
\frac{1}{n},
\end{equation}
but
\begin{equation}\label{0206}
\|v_{n}(t)\|_{\bar{\mathrm{H}}^k}\geq \alpha\ for\ all\ t\geq
t_n:=C_{j,k}+C_1\mathrm{ln}(\alpha n).
\end{equation}Moreover, there exists a positive constant
$\lambda_0$ such that
\begin{equation}\label{0207}
\|\eta_{n}(t)\|_{\bar{\mathrm{H}}^k}\geq\lambda_0\|v_{n}(t)\|_{\bar{\mathrm{H}}^k}\mbox{
and }\|v_{n}(t)\|_{\bar{\mathrm{H}}^k}\rightarrow \infty \mbox{ as }
t\rightarrow \infty.
\end{equation}
\end{thm}

 Theorem \ref{thm:0201} shows globally discontinuous dependence of solutions upon
initial data. The proof of Theorem \ref{thm:0201} is motivated by
\cite{GYTI2} under nontrivial modifications, and its basic idea is
the following. First, noticing that the coefficients of the linearized equations
depend only on the vertical variable $x_3\in (-1,1)$, we seek ``normal mode'' solutions
by taking the horizontal Fourier transform of the equations and assuming that
the solutions grow exponentially in time with the factor
$e^{\lambda(|\xi_1|,|\xi_2|)t}$, where $\xi\in \mathbb{R}^2$ is the
horizontal spatial frequency and $\lambda(|\xi_1|,|\xi_2|)>0$. This
reduces the equations to a system of ordinary differential equations (ODE)
with $\lambda(|\xi_1|,|\xi_2|)>0$ for each $\xi$. Thus, solving the
ODE system by the modified variational method, we show that
$\lambda(|\xi_1|,|\xi_2|)>0$ is a continuous function on some
symmetric open domain $\mathbb{A}\subset \{\xi\in
\mathbb{R}^2~|~|\xi|\in (0,|\xi|_c)\}$ (see Proposition
\ref{pro:0307}), the normal modes with spatial frequency grow in
time, providing a mechanism for the global R-T
instability, where $|\xi|_c=\sqrt{g[\varrho]/\kappa}$ if $\kappa>0$,
otherwise $|\xi|_c=\infty$. In particular, we can restrict $\xi$ in
some symmetric closed sector domain $\mathbb{D}$, such that
$\lambda(|\xi|)$ has a uniform lower bound for any $\xi\in
\mathbb{D}$ (see Remark \ref{rem:0302}). Then, we form a Fourier
synthesis of the normal mode solutions constructed for each spatial frequency
$\xi$ to give solutions of the linearized equations that grow in time, when
measured in $\bar{\mathrm{H}}^k(\Omega)$ for any $k\geq 0$. Finally, we exploit
the property of the trace theorem to show a uniqueness result
of the linearized problem (i.e. Theorem \ref{thm:0401}).
Our results show that in spite of the
uniqueness, the linearized problem is still globally unstable in
$\bar{\mathrm{H}}^k(\Omega)$ for any $k$.

With the linear global instability established, we can show the global instability
of the corresponding nonlinear problem in some sense (i.e. Theorem \ref{thm:0202}).
Recalling that the steady state solution to (\ref{0134}) is
given by $\eta=0$, $v=0$, $q=$constant, we now rewrite the nonlinear
equations (\ref{0132}) in the form of perturbation around the steady
state. Let
\begin{equation*}
\eta=\mbox{Id}+\tilde{\eta}\mbox{ is invertible},\;\;
q=\mbox{constant}+\tilde{q},\;\; v=0+v,\;\; A=I-{G},
\end{equation*}
where
$$  
{G}^T=I-(I+D\tilde{\eta})^{-1}. $$

Then, the evolution equations (\ref{0132}) can be rewritten for
$\tilde{\eta}$, $v$, $\tilde{q}$ as
\begin{equation}\label{0209} \left\{
                              \begin{array}{ll}
\partial_t \tilde{\eta}=v,\\
\varrho\partial_t v_i+(I_{jk}-G_{jk})\partial_k
\tilde{T}_{ij}=\bar{B}_l\bar{B}_m\partial_{lm}^2
\tilde{\eta}_i,\quad i=1,2,3,\\[1mm]
\mathrm{div}v-\mathrm{tr}(G D v)=0,
\end{array}
                            \right.
\end{equation} where tr$(\cdot)$ denotes the matrix trace and
\begin{equation*}
\tilde{T}_{ij}=\tilde{q} I_{ij}-\mu \big((I_{jk}-G_{jk})
\partial_k v_i+(I_{ik}-G_{ik})\partial_k v_j\big).
\end{equation*}
The jump conditions across the interface are
\begin{equation}\label{0211}\begin{array}{l}
\llbracket  v \rrbracket =0,\;\; \llbracket \tilde{T}_{ij}n_j
\rrbracket=g[\varrho]\tilde{\eta}_3 n_i+\bar{B}_l\bar{B}_m\llbracket\partial_l
\tilde{\eta}_i \rrbracket (I_{ij}+\partial_m\tilde{\eta}_j) n_j+\kappa H n_i,\end{array}
\end{equation}
where $n$ and $H$ are respectively given by (\ref{0119}) and
(\ref{0120}) with $\eta=\mathrm{Id}+\tilde{\eta}$. Finally, we
require the boundary conditions
\begin{equation}\label{0212}v_-(t,x',-1)=0,\quad v_+(t,x',1)=0.\end{equation}

We collectively refer to the evolution, jump, and boundary equations
(\ref{0209})--(\ref{0212}) as ``the perturbed problem". To shorten
notation, for $k\geq 0$ we define
\begin{equation*}\label{0213}\|(\tilde{\eta}, v,\tilde{q})(t)\|_{\bar{\mathrm{H}}^k}
=\|\tilde{\eta}(t)\|_{\bar{\mathrm{H}}^k}+
\|v(t)\|_{\bar{\mathrm{H}}^k}
+\|\tilde{q}(t)\|_{\bar{\mathrm{H}}^k}.\end{equation*}
\begin{definition}\label{def:0201}
We say that the perturbed problem has global stability of order $k$
for some $k\geq 3$ if there exist $\delta$, $C_2>0$ and a function
$F:[0,\delta)\to\mathbb{R}^+$ satisfying $F(z)\leq C_2z$ for $z\in
[0,\delta)$, so that the followings hold. For any $T>0$,
$\tilde{\eta}_0$, and $v_0$  satisfying
\begin{equation*}\label{0214}\|(\tilde{\eta}_0,  v_0)\|_{\bar{\mathrm{H}}^k}<\delta,\end{equation*}
there exist $\tilde{\eta}, v,\tilde{q}\in L^\infty (0,T
;\mathrm{\bar{H}}^3(\Omega))$, so that
\begin{enumerate}[\quad\ (1)]
  \item $(\tilde{\eta},  v)(0)=(\tilde{\eta}_0, v_0)$;
  \item $\tilde{\eta}$,  $v$, $\tilde{q}$ solve the perturbed
  problem on $ (0,T)\times\Omega$;
  \item it holds that
\begin{equation*}\label{0215}\sup_{0\leq t\leq T}\|(\tilde{\eta},v,\tilde{q})(t)\|_{\bar{\mathrm{H}}^3}\leq
F(\|(\tilde{\eta}_0,v_0)\|_{\bar{\mathrm{H}}^k});
\end{equation*}
\item $\eta_\pm\in C^2(\bar{\Omega}_\pm)$, respectively, when $\kappa>0$.
\end{enumerate}
\end{definition}
\begin{rem}\label{rem:0201}
In the above definition, in view of (\ref{0119}) and (\ref{0120}),
we have admitted that\newpage
\begin{equation*}n=\frac{\partial_1\eta_\pm\times \partial_2\eta_\pm}{|\partial_1\eta_\pm\times
\partial_2\eta_\pm|}\bigg|_{\{x_3=0\}}
\end{equation*} and
\begin{equation*}H=
\left(\frac{|\partial_1\eta_\pm|^2\partial_2^2\eta_\pm-2(\partial_1\eta_\pm\cdot\partial_2\eta_\pm)\partial_1\partial_2\eta_\pm+
|\partial_2\eta_\pm|^2\partial_1^2
\eta_\pm}{|\partial_1\eta_\pm|^2|\partial_2\eta_\pm|^2-|\partial_1\eta_\pm\cdot\partial_2\eta_\pm|^2}\right)\cdot
n\bigg|_{\{x_3=0\}}
\end{equation*}with $\eta_\pm=\mathrm{Id}+\tilde{\eta}_\pm$.
For convenience in the subsequent proof of Theorem \ref{thm:0202}, we will
still use $\eta$ to denote the both $\eta_-$ and
$\eta_+$ at the interface $\{x_3=0\}$. Similarly, $\tilde{\eta}$ in
(\ref{0211}) also includes the both cases of $\tilde{\eta}_-$ and
$\tilde{\eta}_+$ at $\{x_3=0\}$ except for the notation $\llbracket
\partial_l \tilde{\eta}_i \rrbracket$. We point out that
Theorem \ref{thm:0202} below still holds if $\tilde{\eta}$ only
denotes $\tilde{\eta}_+$ (or $\tilde{\eta}_-$) at $\{x_3=0\}$.
\end{rem}

We can show that the property of global stability of order $k$
cannot hold for any $k\geq 3$, i.e., the following Theorem
\ref{thm:0202}, which will be proved in Section 5.

\begin{thm}\label{thm:0202}
Assume the constant vector $\bar{B}=({B}, 0, 0)$, the perturbed
problem does not have property  of global stability of order $k$ for
any $k\geq 3$.
\end{thm}
The basic idea in the proof of Theorem \ref{thm:0202} is to show, by
utilizing the Lipschitz structure of $F$, that the global stability
of order $k$ would give rise to certain estimates of solutions to
the linearized equations (\ref{0132}) that cannot hold in general
because of Theorem \ref{thm:0201}. We will adapt and modify the
arguments in \cite{GYTI1} to prove Theorem \ref{thm:0202}.
Compared with the perturbed problem in \cite[Theorem 5.2]{GYTI1},
the main difficulty lies in the convergence of the jump
conditions of the rescaled pressure function sequence, because we
could not obtain the strong convergence of the rescaled pressure as in
\cite{GYTI1}. To circumvent such difficulty, similarly to
\cite{JFJSWWWO}, we apply the techniques of integrating by parts to avoid
using the strong convergence of the rescaled pressure. Also,
this idea is used in the proof of the uniqueness of solutions to
the linearized equations (\ref{0134}) in Section 4.

\begin{rem}Theorems 2.1 and 2.2 show that a horizontal magnetic field can not
prevent the linear and nonlinear R-T global instability in the sense
described in Theorems 2.1 and 2.2. From the proof we easily see that
Theorems \ref{thm:0201} and \ref{thm:0202} still hold with
$\bar{B}=(0,B,0)$ in place of $\bar{B}=(B,0,0)$. In addition, based
on the results in \cite{WYC}, we can show that Theorems
\ref{thm:0201} and \ref{thm:0202} still remain valid with
$\bar{B}=(0,0,B)$ in place of $\bar{B}=(B,0,0)$, if one of the
following two conditions holds:
\begin{enumerate}[\quad \ (1)]
  \item $\kappa=0$ and $|B|<g[\varrho]/2$;
  \item $\kappa>0$ and $|B|$ sufficiently small.
\end{enumerate}
\end{rem}

\section{Construction of a growing solution to linearized equations}
\subsection{Growing mode ansatz and the horizontal Fourier transform}

We wish to construct a solution to the linearized equations
(\ref{0134}) that has a growing $H^k$ norm for any $k$. We will
construct such solutions via Fourier synthesis by first constructing
a growing mode for fixed spatial frequency.

 To begin, we assume a
growing mode ansatz (cf. \cite{CSHHS}), i.e.,
\begin{equation*}\label{0301}
v(t,x)=w(x)e^{\lambda t},\;\; q(t,x)=\varpi(x)e^{\lambda t},\;\;
\eta(t,x)=\varsigma(x)e^{\lambda t}
\end{equation*}
for some $\lambda>0$, where $w=(w_1,w_2, w_3)$. Substituting this ansatz into (\ref{0134}),
eliminating $\varsigma$ by using the first equation, we arrive at
the time-independent system for $w$ and $\varpi$:
\begin{equation}\label{0302}
\left\{
                              \begin{array}{ll}
\lambda\varrho w+\nabla \varpi=\mu\Delta w+\lambda^{-1}\bar{B}_l
\bar{B}_m\partial_{lm}^2 w,\\
\mathrm{div}w=0,
\end{array}
                            \right.
\end{equation}
with the corresponding jump conditions
\begin{equation*}\label{0303}
\llbracket  w \rrbracket=0,\;\; \llbracket\varpi I-\mu(Dw+Dw^T)
\rrbracket e_3=\lambda^{-1}(g[\varrho]w_3+\kappa\Delta_{x'}w_3 )e_3
+\lambda^{-1}\bar{B}_3\bar{B}_l\llbracket  \partial_l w \rrbracket,
\end{equation*}
and the boundary conditions
\begin{equation*}\label{0304}
w(t,x',-1)=0,\quad w(t,x',1)=0.
\end{equation*}

 We take the horizontal Fourier transform of $w_1$, $w_2$, $w_3$ in (\ref{0302}),
which we denote by either $\hat{\cdot}$ or $\mathcal{F}$, and fix
a spatial frequency $\xi=(\xi_1, \xi_2)\in \mathbb{R}^2$. Define the
new unknowns $\varphi(x_3)=i\hat{w}_1(\xi,x_3)$,
$\theta(x_3)=i\hat{w}_2(\xi,x_3)$, $\psi(x_3)=\hat{w}_3(\xi,x_3)$
and $\pi(x_3)=\hat{\varpi}_3(\xi,x_3)$, so that
\begin{equation*}\label{0305}
\mathcal{F}(\mathrm{div}w)=\xi_1\varphi+\xi_2\theta+\psi',
\end{equation*}
where $'=d/dx_3$. Since we consider the case
$\bar{B}=(\bar{B}_1,\bar{B}_2,\bar{B}_3)\equiv (B,0,0)$, then for
$\varphi$, $\theta$, $\psi$ and $\lambda=\lambda(\xi)$ we arrive at
the following system of ODEs.
\begin{equation}\label{0306} \left\{
                              \begin{array}{ll}
\lambda^2\varrho \varphi-\lambda\xi_1\pi
+ \lambda\mu(|\xi|^2\varphi-\varphi'')+B^2\xi_1^2\varphi=0,\\[1mm]
\lambda^2\varrho \theta-\lambda\xi_2\pi+\lambda\mu(|\xi|^2\theta-\theta'')+B^2\xi_1^2\theta=0,\\[1mm]
\lambda^2\varrho \psi+\lambda\pi'+\lambda\mu(|\xi|^2\psi-\psi'')+B^2\xi_1^2\psi=0,\\[1mm]
\xi_1\varphi+\xi_2\theta+\psi'=0,
\end{array}
                            \right.
\end{equation}
along with the jump conditions
\begin{equation}\label{0307} \left\{
                              \begin{array}{ll}
\llbracket  \varphi \rrbracket=\llbracket
\theta\rrbracket=\llbracket\psi \rrbracket=0,\\[1mm]
\llbracket\lambda\mu (\xi_1\psi-\varphi') \rrbracket=\llbracket\lambda\mu
(\xi_2\psi-\theta') \rrbracket=0, \\[1mm]
 \llbracket -2\lambda\mu \psi'+\lambda\pi \rrbracket=(g[\varrho]-k|\xi|^2)\psi,
\end{array}
                            \right.
\end{equation}
and the boundary conditions
\begin{equation}\label{0308}
\varphi(-1)=\varphi(1)=\theta(-1)=\theta(1)=\psi(-1)=\psi(1)=0.
\end{equation}

Eliminating $\pi$ from the third equation in (\ref{0306}), we obtain
the following ODE for $\psi$
\begin{equation}\label{0309}
-\lambda^2\rho(|\xi|^2\psi-\psi'')=\lambda\mu(|\xi|^4\psi-2|\xi|^2\psi''+\psi'''')+B^2\xi_1^2(|\xi|^2\psi-\psi'')
\end{equation}
along with the jump conditions
\begin{eqnarray}
&&\llbracket\psi \rrbracket=\llbracket\psi' \rrbracket=0,\label{0310}\\
&&\llbracket \lambda\mu (|\xi|^2\psi+\psi'')
\rrbracket=0,\label{0311}
\\&& \llbracket
\lambda\mu(\psi'''-3|\xi|^2\psi')\rrbracket=\llbracket
\lambda^2\varrho\psi'\rrbracket+(g[\rho]-k|\xi|^2)|\xi|^2\psi,\label{0312}
\end{eqnarray}
and the boundary conditions
\begin{equation}\label{0313}
\psi(-1)=\psi(1)=\psi'(-1)=\psi'(1)=0.
\end{equation}

Before constructing a growing solution to ODE of
(\ref{0309})--(\ref{0313}) in next subsection, we should introduce
so-called admissible set for a growing solution.

First we define the maximal frequency by
  \begin{equation}\label{0314}
|\xi|_c:=\left\{
  \begin{array}{ll}
     \sqrt{{g[\varrho]}/{\kappa}}, & \hbox{ if }\kappa>0,  \\[0.3em]
    +\infty, & \hbox{ if }\kappa=0,
  \end{array}
\right.
 \end{equation}
and the generalized critical magnetic number (depending on the frequency) by
\begin{equation}\label{0315}
B_c(|\xi|):=\sqrt{\sup_{\psi \in H^1_0(-1,1)\atop \psi\not\equiv 0}
\frac{(g[\varrho]-\kappa|\xi|^2)\psi^2(0)}{\int_{-1}^1|\psi'|^2\mathrm{d}y}}\quad\;\;
\mbox{ for }|\xi|<|\xi|_c,
\end{equation}
and the generalized critical frequency function by
\begin{equation}\label{0316}
S(|\xi_1|,|\xi_2|)=\sqrt{\sup_{\psi\in
H^1_0{(-1,1)}\atop\psi\not\equiv 0}
\frac{(g[\varrho]-\kappa|\xi|^2)\psi^2(0)-f^2(|\xi_1|,|\xi_2|)
\int_{-1}^1|\psi'|^2\mathrm{d}y}{f^2(|\xi_1|,|\xi_2|)\int_{-1}^1|\psi|^2\mathrm{d}y}}
\end{equation}
for $|\xi|<|\xi|_c$ and $0<f(\xi)<B_c(|\xi|)$, where
$f(|\xi_1|,|\xi_2|)=|\xi_1B|/|\xi|$ and $S(|\xi_1|,|\xi_2|):=\infty$
when $\xi_1 B=0$.

By \cite[Lemma 3.2]{WYC}, we see that the supremums in both (\ref{0315}) and
 (\ref{0316}) are achieved for each $\xi$ and $B$ with
$\xi_1 B\neq 0$. So, (\ref{0315}) and (\ref{0316}) make sense.
Moreover $S(|\xi_1|,|\xi_2|)\rightarrow \infty$ as
$f(|\xi_1|,|\xi_2|)\rightarrow 0$.

 Now, we define an admissible set for a growing solution
\begin{equation}\label{0319}
\mathbb{A}=\big\{\xi\in \mathbb{R}^2~|~|\xi|<|\xi|_c,\
0<f(|\xi_1|,|\xi_2|)<B_c(|\xi|),\ 0<|\xi|<S(|\xi_1|,|\xi_2|)\big\}.
\end{equation}
By construction, the admissible set possesses the following properties,
which will be useful in Subsections 3.2--3.4.
\begin{pro}\label{pro:0301} Assume $\mathbb{A}$ is defined by (\ref{0319}), then
\begin{enumerate}[\quad \ (1)]
  \item the set $\mathbb{A}$ is a open set in $\mathbb{R}^2$;
  \item the set $\mathbb{A}$ is symmetric on $x$-axis and $y$-axis
in $\mathbb{R}^2$, respectively;
  \item  there exist three positive constants $d_1<d_2$
and $d_3$, such that
\begin{equation}\label{0320} \mathbb{D}:=\{\xi\in
\mathbb{R}^2~|~d_1\leq |\xi|\leq d_2,\ |\xi_1|\leq d_3\}\subset
\mathbb{A}.
\end{equation}
\end{enumerate}
\end{pro}
\begin{pf}The proof is very elementary and we omit it here.
\end{pf}

 In the next subsection, we will show that for any
fixed $B$ and any $\xi\in \mathbb{A}$, there is a nontrivial
solution $\psi$ with $\lambda>0$ to the problem
(\ref{0309})--(\ref{0313}). The next proposition shows that
$B_c(|\xi|)$ can be in fact represented by $g$, $[\varrho]$,
$\kappa$ and $|\xi|$.
\begin{pro}\label{pro:0302} Let $B_c(|\xi|)$ be defined by (\ref{0315}). Then,
$B_c(|\xi|)=\sqrt{(g[\varrho]-\kappa|\xi|^2)/2}$.
\end{pro}
\begin{pf} By virtue of (\ref{0315}), it suffices to prove that
\begin{equation}\label{0321}
{\sup_{\psi \in
H^1_0(-1,1)}\frac{\psi^2(0)}{\int_{-1}^1|\psi'|^2\mathrm{d}y}}=\frac{1}{2}.
\end{equation}

 For any $\psi\in H^1_0(-1,1)$, we have
 $\psi(0)=\int_{-1}^0\psi'\mathrm{d}y=\int_{1}^0\psi'\mathrm{d}y$, which gives
$$  |\psi(0)|\leq \frac{1}{2}\int_{-1}^1|\psi'|\mathrm{d}y. $$
Hence, by the H\"{o}lder inequality, we see that
\begin{equation}\label{0323}
\frac{\psi^2(0)}{\int_{-1}^1|\psi'|^2\mathrm{d}y}\leq \frac{1}{4}
\frac{\left(\int_{-1}^1|\psi'|\mathrm{d}y\right)^2}{\int_{-1}^1|\psi'|^2\mathrm{d}y}
\leq \frac{1}{2}.
\end{equation}
On the other hand, it is easy to check that the function
$$ \psi_s=\left\{
                              \begin{array}{ll}
1+x, & x\in (-1,0],\\
1-x, & x\in (0,1)
\end{array}  \right. $$
satisfies $\psi_s\in H^1_0(-1,1)$. Furthermore, a simple computation results in
$$   \frac{\psi^2_s(0)}{\int_{-1}^1|\psi'_s|^2\mathrm{d}y}=\frac{1}{2}, $$
which combined with (\ref{0323}) implies (\ref{0321}).
\hfill $\Box$
\end{pf}
\begin{rem}\label{rem:03011111}
We mention that the definition of critical number and critical
frequency (both are independent of frequency $\xi$) was introduced
by Wang \cite{WYC}, where he defined the  critical number by
\begin{equation*}
|B|_c:=\sqrt{\sup_{\psi \in
H^1_0(-1,1)}\frac{g[\varrho]\psi^2(0)}{\int_{-1}^1|\psi'|^2\mathrm{d}y}},
\end{equation*}
which is equal to $\sqrt{g[\varrho]/2}$ by virtue of Proposition \ref{pro:0302}.
\end{rem}
\subsection{Construction of a growing solution to ODE}
 Throughout this subsection we assume that $\xi\in \mathbb{A}$. Similarly to
\cite{GYTI2,WYC}, we can apply the variational method to construct
a growing solution of (\ref{0309})--(\ref{0313}) for given $\xi\in
\mathbb{A}$. For the reader's convenience, we sketch the procedure of
the construction.

First of all, in order to apply the variational method, we
remove the linear dependence on $\lambda$ in (\ref{0309})--(\ref{0312}) by defining the modified
viscosities $\bar{\mu}=s\mu$, where $s>0$ is an arbitrary parameter.
Thus we obtain a family ($s>0$) of modified problems
\begin{equation}\label{0326}
-\lambda^2\rho(|\xi|^2\psi-\psi'')=s\mu
(|\xi|^4\psi-2|\xi|^2\psi''+\psi'''')+B^2\xi_1^2(|\xi|^2\psi-\psi'')
\end{equation}
along with the jump conditions
\begin{eqnarray}
&&\llbracket\psi \rrbracket=\llbracket\psi' \rrbracket=0,\label{0327}\\
&&\llbracket s\mu (|\xi|^2\psi+\psi'') \rrbracket=0,\label{0328}
\\&& \llbracket
s \mu(\psi'''-3|\xi|^2\psi')\rrbracket=\llbracket
\lambda^2\varrho\psi'\rrbracket+(g[\rho]-\kappa|\xi|^2)|\xi|^2\psi,\label{0329}
\end{eqnarray}
and the boundary conditions
\begin{equation}\label{0330}
\psi(-1)=\psi(1)=\psi'(-1)=\psi'(1)=0.
\end{equation}
Notice that for any fixed $s>0$ and $\xi$,
(\ref{0326})--(\ref{0330}) is a standard eigenvalue problem for
$-\lambda^2$. This allows us to use the variational method to
construct solutions. We define the energies by
\begin{equation}\label{0331}\begin{aligned}
E(\psi)=&\frac{1}{2}\int_{-1}^1s\mu(4|\xi|^2|\psi'|^2+|\xi|^2\psi+\psi''|^2)+B^2\xi_1^2(|\psi'|^2+
|\xi|^2\psi^2)\mathrm{d}x_3\\
&- \frac{1}{2}|\xi|^2(g[\rho]-\kappa|\xi|^2)\psi^2(0),
\end{aligned}\end{equation}
and
\begin{equation}\label{0332}
J(\psi)=\frac{1}{2}\int_{-1}^1\rho(|\xi|^2\psi^2+|\psi'|^2)\mathrm{d}x_3\end{equation}
which are well-defined on the space $H_0^2(-1,1)$. We define the admissible
set for the energy (\ref{0331})
\begin{equation*}\label{}
\mathcal{A}=\{\psi\in H^2_0(-1,1)~|~J(\psi)=1\}.  \end{equation*}
Thus we can find the smallest $-\lambda^2$ by minimizing
\begin{equation}\label{0334}
-\lambda^2(|\xi_1|,|\xi_2|)=\alpha(|\xi_1|,|\xi_2|):=\inf_{\psi\in
\mathcal{A}}E(\psi).\end{equation} In fact, we can show that a
minimizer of (\ref{0334}) exists and satisfies the
Euler-Lagrange equations equivalent to (\ref{0326})--(\ref{0330}).

\begin{pro}\label{pro:0303}
 For any fixed $\xi\in \mathbb{A}$, $E$ achieves its infinimum on
$\mathcal{A}$. In addition, let $\psi$ be a minimizer and
$-\lambda^2:=E(\psi)$, then the pair ($\psi$, $\lambda^2$) satisfies
(\ref{0326}) along with the jump and boundary conditions (\ref{0327})--(\ref{0330}).
Moreover, $\psi$ is smooth when restricted to ($-1,0$) or ($0,1$).
\end{pro}
\begin{pf} We can follow the same proof procedure as in \cite[Proposition 3.1]{WYC} to
obtain Proposition \ref{pro:0301}. Hence, we omit the details of the
proof here.\hfill $\Box$
\end{pf}

 In order to emphasize the
dependence on $s\in (0,\infty)$, we will sometimes write
\begin{equation*}
\lambda(|\xi_1|,|\xi_2|)=\lambda(|\xi_1|,|\xi_2|,s)\end{equation*}
and
\begin{equation*}
\alpha(|\xi_1|,|\xi_2|)=\alpha(s).\end{equation*} Next, we want to
prove that there is a fixed point such that $\lambda=s$. To this
end, we first give some properties of $\alpha(s)$ as a function of
$s> 0$ for each fixed $\xi\in \mathbb{A}$.

\begin{pro}\label{pro:0304} Let
$\alpha(s):(0,\infty)\to\mathbb{R}$ be given by
(\ref{0334}), then the following conclusions hold:
\begin{enumerate}[\quad \ (1)]
  \item  $\alpha(s)\in C_{\mathrm{loc}}^{0,1}(0,\infty)$ is strictly
  increasing;
  \item there exist two positive constants $c_1:=c_1(g,\varrho_\pm)$ and $c_2:=c_2(|\xi|,\mu_\pm)$,
   such that
  \begin{equation*}\alpha(s)\geq -c_1|\xi|+sc_2;
\end{equation*}
 \item for each fixed $\xi\in\mathbb{A}$, there exist two positive constants $c_3:
 =c_3(|B|,|\xi|, g,\varrho_\pm,\kappa)$ and $c_4:=c_4(|B|,|\xi|, g,\varrho_\pm,\mu_\pm,\kappa)$,
 such that
  \begin{equation}\label{0336}\alpha(s)\leq -c_3+sc_4.\end{equation}
\end{enumerate}
\end{pro}
\begin{pf} The first two assertions can be shown in the same way as in \cite[Lemma 3.5]{WYC} and
\cite[Propostition 3.6]{GYTI2}. It remains to prove the third assertion to complete the proof.

First, observe that the energy $E(\psi)$ can be decomposed in the form:
\begin{equation*}E(\psi)=|\xi|^2E_0(\psi)+sE_1(\psi),
\end{equation*}
where
\begin{eqnarray*}
&&E_0(\psi)=\frac{1}{2}\int_{-1}^1
f^2(|\xi_1|,|\xi_2|)(|\psi'|^2+|\xi|^2\psi^2)\mathrm{d}x_3-\frac{1}{2}
\left(g[\varrho]-\kappa|\xi|^2\right)\psi^2(0),\label{}\\
&&E_1(\psi)=\frac{1}{2}\int_{-1}^1\mu
(4|\xi|^2|\psi'|^2+||\xi|^2\psi+\psi''|^2)\mathrm{d}x_3,\\
&& f(|\xi_1|,|\xi_2|)=\frac{\left|\xi_1 B\right|}{|\xi|}.
\end{eqnarray*}
Then, since $\xi\in \mathbb{A}$ by virtue of the definitions
(\ref{0316}) and (\ref{0319}), there exists $\bar{\psi}\in H_0^1$ such that
\begin{equation}\label{0341}
\int_{-1}^1f^2(|\xi_1|,|\xi_2|)(|\xi|^2\bar{\psi}^2+|\bar{\psi}'|^2)\mathrm{d}x_3<
\left(g[\varrho]-\kappa|\xi|^2\right)\psi^2(0).
\end{equation}
On the other hand, $C_0^\infty(-1,1)$ is dense in $H_0^1(-1,1)$,
then there exists a function sequence $\bar{\psi}_n\in C_0^\infty(-1,1)$, so that
\begin{equation}\label{0342}\bar{\psi}_n\rightarrow\bar{\psi}\mbox{ stronly in
}H_0^1(-1,1),
\end{equation}
which, together with the compact embedding theorem,  yields
\begin{equation}\label{0343}\bar{\psi}_{n}\rightarrow\bar{\psi}\mbox{ stronly in
}C^0(-1,1).
\end{equation}
Putting (\ref{0341})--(\ref{0343}) together, we see that there exist a subsequence
$\bar{\psi}_{n_0}\in H_0^2(-1,1)$, such that $\bar{\psi}_{n_0}\not\equiv 0$ and
\begin{equation}\label{0344}
\int_{-1}^1f^2(|\xi_1|,|\xi_2|)(|\xi|^2|\bar{\psi}_{n_0}|^2+|\bar{\psi}_{n_0}'|^2)\mathrm{d}x_3<
\left(g[\varrho]-\kappa|\xi|^2\right)\bar{\psi}_{n_0}^2(0),
\end{equation}
which implies
\begin{equation*} E_0(\bar{\psi}_{n_0})<0. \end{equation*}
Thus, we have
\begin{equation}\begin{aligned} \label{0345}
\alpha(s)=&\inf_{\psi\in \mathcal{A}}E(\psi)=\inf_{\psi\in
H_0^2(-1,1)\atop \psi\equiv\!\!\!\! /\
0}\frac{E(\psi)}{J(\psi)}\\
\leq &
\frac{E(\bar{\psi}_{n_0})}{J(\bar{\psi}_{n_0})}=|\xi|^2\frac{E_0(\bar{\psi}_{n_0})}{J(\bar{\psi}_{n_0})}
+s\frac{E_1(\bar{\psi}_{n_0})}{J(\bar{\psi}_{n_0})}:= -c_3+sc_4
\end{aligned}\end{equation} for two positive constants $c_3:
 =c_3(|B|,|\xi|, g,\varrho_\pm,\kappa)$ and $c_4:=c_4(|B|,|\xi|, g,\varrho_\pm,\mu_\pm,\kappa)$.
 \hfill $\Box$  \end{pf}

Given $\xi\in \mathbb{A}$, by virtue of (\ref{0336}), there exists a
$s_0>0$ depending on the quantities $|B|,|\xi|, g,\varrho_\pm$,
$\mu_\pm$ and $\kappa$, so that for $s\leq s_0$, $\alpha(s)<0$.
Therefore, we can define the open set
\begin{equation*}\mathcal{S}=\alpha^{-1}(-\infty,0)\subset (0,\infty).
\end{equation*}
Note that $\mathcal{S}$ is non-empty and this allows us to define
$\lambda(s)=\sqrt{-\alpha(s)}$ for $s\in \mathcal{S}$. Thus, as
a result of Proposition \ref{pro:0301}, we have the following
existence result for the modified problem (\ref{0326})--(\ref{0330}).
\begin{pro}\label{pro:0305}
For each $\xi\in \mathbb{A}$ and each $s\in \mathcal{S}$, there is a
solution $\psi=\psi(|\xi_1|,|\xi_2|, x_3)$ with
$\lambda=\lambda(|\xi_1|,|\xi_2|, s)>0$ to the problem (\ref{0326})
with the jump and boundary conditions (\ref{0327})--(\ref{0330}).
Moreover, $\psi$ is smooth when restricted to ($-1,0$) or ($0,1$)
with $\psi(|\xi_1|,|\xi_2|,0)\neq 0$.
\end{pro}

 Finally, we can use Proposition \ref{pro:0304} to make a
fixed-point argument to find a $s\in \mathcal{S}$ satisfying
$s=\lambda(|\xi_1|,|\xi_2|, s)$ and to construct solutions to the
original problem (\ref{0306})--(\ref{0308}).
\begin{pro}\label{pro:0306} Let $\xi\in\mathbb{A}$,
then there exists a unique $s\in \mathcal{S}$, so that
$\lambda(|\xi_1|,|\xi_2|, s)=\sqrt{-\alpha(s)}>0$ and $s=\lambda(|\xi_1|,|\xi_2|, s)$.
\end{pro}
\begin{pf}
We refer to \cite[Lemma 3.7]{WYC} for a proof.\hfill $\Box$
\end{pf}

 Consequently, in view of Propositions \ref{pro:0305} and \ref{pro:0306},
 we conclude the following existence result on the problem (\ref{0306})--(\ref{0308}).
\begin{thm}\label{thm:0301}
For each $\xi\in \mathbb{A}$, there exist
$\psi=\psi(|\xi_1|,|\xi_2|,x_3)$ and $\lambda(|\xi_1|,|\xi_2|)>0$
satisfying
 (\ref{0306})--(\ref{0308}). Moreover, $\psi$ is smooth when restricted to ($-1,0$) or
  ($0,1$) with $\psi(|\xi_1|,|\xi_2|,0)\neq 0$.
\end{thm}

Next, we show some properties of the solutions established in
Theorem \ref{thm:0301} in terms of
$\lambda(\xi):=\lambda(|\xi_1|,|\xi_2|)$. The first property is
given in the following proposition which shows that $\lambda$ is a
bounded, continuous function of $\xi\in \mathbb{A}$. To this end, we
shall apply the following Ehrling-Nirenberg-Gagliardo
interpolation inequality, the proof of which can be found in
\cite[Chapter 5]{ARAJJFF}.
\begin{lem}\label{lem:0301}
Let $\Omega$ be a domain in $\mathbb{R}^n$ satisfying the cone
condition. For each $\varepsilon_0>0$ there is a constant $K$
depending on $n$, $m$, $p$ and $\varepsilon_0$,
such that for $0<\varepsilon\leq \varepsilon_0$, $0\leq j\leq m$ and
$u\in W^{m,p}(\Omega)$,
\begin{equation*}
\sum_{|\alpha|=j}\int_\Omega|D^\alpha u(x)|^p\mathrm{d}x\leq
K\Big(\varepsilon\sum_{|\alpha|=m}\int_\Omega|D^\alpha
u(x)|^p\mathrm{d}x +\varepsilon^{-j/(m-j)}\int_\Omega
|u|^p\mathrm{d}x\Big).
\end{equation*}
\end{lem}

\begin{pro}\label{pro:0307} The function
$\lambda :\mathbb{A}\rightarrow (0,\infty)$ is continuous and
satisfies
\begin{equation}\label{0348}
\Lambda:=\sup_{\xi\in \mathbb{A}}\lambda(\xi)\leq
\frac{g[\varrho]}{4\mu_-}.\end{equation} Moreover,
\begin{equation}\label{0349}
\lim_{\xi\in\mathbb{A},\xi\to 0}\lambda(\xi)=0,\end{equation}
and if $\kappa>0$, then also
\begin{equation}\label{0350} \lim_{ \xi\in \mathbb{A},\xi\to\xi_c}
\lambda(\xi)=0\;\; \mbox{ for any }\,\xi_c\in\bar{\mathbb{A}}\,
\mbox{ with }\,|\xi_c|=|\xi|_c,
\end{equation}
where $\bar{\mathbb{A}}$ denotes the closure of
$\mathbb{A}$ and $|\xi|_c$ is defined in (\ref{0314}).
\end{pro}
\begin{pf} The proposition can be shown by a proof procedure similar to that
used in \cite[Proposition 3.9]{GYTI2} but with necessary nontrivial modifications
in arguments. For the reader's convenience, we give the proof in details here.

(i) We start with proving the boundedness of $\lambda$. By (\ref{0334}), we find that
$$ 
\begin{aligned}
-\lambda^2(\xi)=&\frac{1}{2}\int_{-1}^1\lambda(\xi)\mu(4|\xi|^2|
\psi'_{\xi}|^2+||\xi|^2\psi_{\xi}+\psi''_{\xi}|^2)\\
&+|B|^2\xi_1^2(|\psi'_{\xi}|^2+|\xi|^2\psi^2_{\xi})\mathrm{d}x_3
-\frac{1}{2}|\xi|^2(g[\rho]-\kappa|\xi|^2)\psi^2_{\xi}(0),
\end{aligned} $$ 
which yields
\begin{equation}\label{0352}
2\mu_-|\xi|^2\lambda(\xi)\int_{-1}^1|\psi'_{\xi}|^2
\mathrm{d}x_3\leq\frac{1}{2}|\xi|^2g[\rho]\psi_{\xi}^2(0).
\end{equation}
Using the H\"{o}lder inequality, we can bound
$$  \psi_{\xi}^2(0)=\left|\int_0^1\psi'_{\xi}\mathrm{d}x_3\right|^2 \leq
\int_{0}^1|\psi'_{\xi}|^2\mathrm{d}x_3. $$
Substitution of the above inequality into (\ref{0352}) gives then
\begin{equation}\label{0354}
|\xi|^2\left(2\mu_-\lambda(|\xi|)-\frac{1}{2}g[\varrho]\right)
\int_{-1}^1|\psi'_{\xi}|^2\mathrm{d}x_3\leq 0.
\end{equation}
Consequently, (\ref{0354}) implies (\ref{0348}), since
$\|\psi'_{\xi}\|_{L^2(-1,1)}>0$.

(ii) We now turn to the proof of the continuity claim. Since
$\lambda=\sqrt{-\alpha}$, it suffices to show the continuity of
$\alpha(\xi):=\alpha(|\xi_1|,|\xi_2|)$. By virtue of Theorem \ref{pro:0301},
for each $\xi\in\mathbb{A}$ there exists a function $\psi_{\xi}\in \mathcal{A}$
satisfying (\ref{0306})--(\ref{0308}), so that
$\alpha(\xi)=E(\psi_{\xi})$. Furthermore, $\psi$ is smooth when
restricted to ($-1,0$) or ($0,1$).
We have that $\alpha(\xi)<0$, which, when combined with
(\ref{0331}), yields the estimate
\begin{equation}\label{0355}\begin{aligned}
\frac{1}{2}\int_{-1}^1s\mu||\xi|^2\psi_\xi+\psi''_\xi|^2\mathrm{d}x_3-
\frac{1}{2}|\xi|^2(g[\rho]-\kappa|\xi|^2)\psi_\xi^2(0)\leq\alpha (\xi )<0.
\end{aligned}\end{equation}
On the other hand,
\begin{equation}\label{0356}|\xi|\psi_\xi^2(0)\leq
2\Big(\int_0^1|\xi|^2\psi^2_\xi\mathrm{d}x_3\Big)^{1/2}
\Big(\int_0^1 |\psi'_\xi|^2\mathrm{d}x_3\Big)^{1/2}\leq
\frac{4}{\varrho_-}\end{equation}
because of $\psi_{\xi}\in\mathcal{A}$.
Thus, plugging (\ref{0356}) into (\ref{0355}), we get
\begin{equation}\label{0357}
\frac{1}{2}\int_{-1}^1s\mu||\xi|^2\psi_\xi+\psi''_\xi|^2\mathrm{d}x_3<
{2}|\xi|(g[\rho]-\kappa|\xi|^2)/{\varrho_-}.
\end{equation}

 Now suppose $\xi_n\in
\mathbb{A}$ is a sequence so that $\xi_n\rightarrow \xi\in
\mathbb{A}$. Since $\mathbb{A}$ is a open set and
$\xi\neq(0,0)$, when $n$ is sufficiently large, there exists
a sufficiently small, bounded and open sector domain
$\mathbb{S}\subset \mathbb{A}$ satisfying the following three
conditions:
\begin{itemize}
  \item there exists a $n_0>0$ such that $\xi, \xi_n\in\mathbb{S}$ for any $n>n_0$;
  \item $(0,0)\in \!\!\!\!\!/\ \bar{\mathbb{S}}$, where $\bar{\mathbb{S}}$
  denotes the closure of $\mathbb{S}$;
  \item there is a $\bar{\xi}\in \bar{\mathbb{S}}$ such that
\begin{equation}\label{0359}
\frac{|\xi_1|}{|\xi|}\leq \frac{|\bar{\xi}_1|}{|\bar{\xi}|},\;\;
|\xi|\leq {|\bar{\xi}|}\;\;\mbox{ for any }\,\xi\in {\mathbb{S}}.
\end{equation}
\end{itemize}

 In order to make use of (\ref{0357}) we have to show that $s(\xi)$ is
bounded uniformly from below for $n>n_0$. Employing arguments similar to those used in
the derivation of (\ref{0344}), we find that there is a $\bar{\psi}\in H_0^2(-1,1)$, such that
$\bar{\psi}\not\equiv 0$ and
$$
\bar{E}_0(\bar{\psi}):=\frac{1}{2}\int_{-1}^1\frac{|\bar{\xi}_1B|^2}{|\bar{\xi}|^2}
(|\bar{\psi}'|^2+|\bar{\xi}|^2\bar{\psi}^2) \mathrm{d}x_3-\frac{1}{2}
(g[\varrho]-\kappa|\bar{\xi}|^2)\bar{\psi}^2(0)<0, $$
which, together with (\ref{0359}),  implies
\begin{equation}\label{0361}
\frac{1}{2}\int_{-1}^1\frac{|{\xi}_1B|^2}{|{\xi}|^2}(|{\xi}|^2\bar{\psi}^2+|\bar{\psi}'|^2)
\mathrm{d}x_3-\frac{1}{2}
(g[\varrho]-\kappa|\xi|^2)\bar{\psi}^2(0)<\bar{E}_0(\bar{\psi})\;\;
\mbox{ for any }\xi\in {\mathbb{S}}.
\end{equation}

Let $|\xi|_0=\inf_{n>n_0}|\xi|$, then $|\xi|_0>0$. Denoting
\begin{eqnarray*}&&
\bar{J}(\bar{\psi})=\frac{1}{2}\int_{-1}^1\varrho(|\bar{\xi}|^2\bar{\psi}^2+|\bar{\psi}'|^2)\mathrm{d}x_3,
\;\;
{J}_c(\bar{\psi})=\frac{1}{2}\int_{-1}^1\varrho(|{\xi}|^2_c\bar{\psi}^2+|\bar{\psi}'|^2)\mathrm{d}x_3,\\
&&\bar{E}_1(\psi)=\frac{1}{2}\int_{-1}^1\mu
(4|\bar{\xi}|^2|\bar{\psi}'|^2+2|\bar{\xi}|^4\psi^2+2|\bar{\psi}''|^2)\mathrm{d}x_3,
\end{eqnarray*}
and using (\ref{0361}), we argue, similarly to (\ref{0345}), to deduce that
$$  
  \alpha(s)\leq -c_1+sc_2\quad\mbox{for any }\xi\in {\mathbb{S}}, $$
where the two constants
$c_1=c_1(|\xi|_0,\bar{E}_0(\bar{\psi}), \bar{J}(\bar{\psi}))$ and
$c_2=c_2(\bar{E}_1(\bar{\psi}), {J}_c(\bar{\psi}))$ are independent
of $s$. Keeping in mind that $ -\alpha(\xi_n)=\lambda^2(\xi_n)=s^2(\xi_n)$, one gets
\begin{equation*}\label{0363}
0\leq s^2(\xi_n)+c_2s(\xi_n)-c_1\;\;\mbox{ for any }n>n_0.
\end{equation*}
 On the other hand, the fact that $\psi\in \mathcal{A}$ shows that $\psi_{\xi_n}$
are uniformly bounded in $H^1(-1,1)$ for any $n>n_0$. Consequently, by virtue of
(\ref{0359}) and (\ref{0327}), the estimate (\ref{0357}) implies
the uniform boundedness of $\psi_{\xi_n}$ in $H^2_0(-1,1)$ for any $n>n_0$.

Plugging the above boundedness facts on $\psi_{\xi_n}$ into the ODE (\ref{0309}) in the intervals
$(-1,0)$ and $(0,1)$ respectively, we find that $\psi_{\xi_n}''''$ are also
uniformly bounded in $L^2(-1,0)$ and $L^2(0,1)$ for any $n>n_0$.
Thus, using Lemma \ref{lem:0301}, we infer that $\psi_{\xi_n}$ are
uniformly bounded in $H^4(-1,0)$ and $H^4(0,1)$ for any $n>n_0$. So,
up to extraction of a subsequence we conclude that
$$ \psi_{\xi_n}\rightarrow \psi\;\;\mbox{ strongly in }\; H^2_0(-1,1),\;
H^3(-1,0)\;\mbox{ and }\; H^3(0,1), $$
which yields
\begin{equation} \label{0364}
\alpha(\xi_n)=E(\psi_{\xi_n})\rightarrow E(\psi_{\xi})=\alpha(\xi).
\end{equation}
Since (\ref{0364}) must hold for any such extracted
subsequence, one deduces that $\alpha(\xi_n)\rightarrow \alpha(\xi)$
for the original sequence $\xi_n$ as well, and hence $\alpha(\xi)$
is continuous.

Finally, we derive the limits as $\xi\rightarrow 0$ when $\kappa\geq
0$, and $\xi\rightarrow \xi_c\in \bar{\mathbb{A}}$ with
$|\xi_c|=|\xi|_c$ when $\kappa>0$, where we restrict the variables
$\xi$ in $\mathbb{A}$. By virtue of (\ref{0334}) and (\ref{0356}), $0\leq
\lambda^2(\xi)\leq 2g[\varrho]|\xi|/\varrho_-$, which gives (\ref{0349}).
On the other hand, when $\kappa>0$ we may utilize (\ref{0334}) to find that
\begin{equation*}
\lambda^2(\xi)\leq {(g[\varrho]-\kappa
|\xi|^2)|\xi|^2}\psi_{\xi}^2(0)/{2}\leq
{2|\xi|(g[\varrho]-\kappa|\xi|^2)}/{\varrho_-},
\end{equation*}which results in (\ref{0350}) when $\kappa>0$.
 \hfill   $\Box$
\end{pf}
\begin{rem}\label{rem:0302111111}
The stabilizing effect of viscosity and surface tension
is evident from the above calculations. As shown in \cite{DRJFJS},
without viscosity or surface tension, there exists a domain
$\mathbb{D}$, such that $\lambda(\xi)\rightarrow \infty$ for $\xi\in
\mathbb{D}$ as $|\xi|\rightarrow \infty$. With viscosity but no
surface tension, by virtue of the definition (\ref{0319}), there are only
partial spatial frequencies $\xi$ (including all
$\xi$ with $\xi_1=0$ in particular) which remain unstable, but the growth of
$\lambda(\xi)$ is bounded. With both viscosity and surface tension, only
those spatial frequencies $\xi$ belonging to the bounded domain
$\mathbb{A}$ are unstable, and $\lambda(\xi)$ remains bounded. In
addition, in the construction of the normal mode solution to the
linearized system, the horizontal magnetic field can enhance the
maximal growth rate $\lambda(\xi)$. In particular, when $\xi_1=0$,
the growth rate reduces to the one for the corresponding equations
of incompressible viscous fluids.
\end{rem}
\begin{rem}\label{rem:0302} From Proposition \ref{pro:0307}, we immediately infer that
\begin{equation}\label{0367}
\lambda_0:=\inf_{\xi\in \mathbb{D}}\lambda(\xi)>0,
\end{equation}where the closed domain $\mathbb{D}$ is defined by (\ref{0320}).
\end{rem}

\subsection{Construction of a solution to the system (\ref{0306})--(\ref{0308})}
 A solution to (\ref{0309})--(\ref{0313}) gives rise to a solution of the system
(\ref{0306})--(\ref{0308}) for the growing mode velocity $v$, as well.
\begin{thm}\label{thm:0302}
For each $\xi\in \mathbb{A}$, there exists a solution
$(\tilde{\varphi},\tilde{\theta},\tilde{\psi},\tilde{\pi})=
\big(\tilde{\varphi}(\xi, x_3),\tilde{\theta}(\xi, x_3)$,
$\tilde{\psi}(\xi,x_3)$, $\tilde{\pi}(\xi ,x_3)\big)$ with
$\lambda =\lambda(|\xi_1|,|\xi_2|)>0$ to (\ref{0306})--(\ref{0308}), and the
solution is smooth when restricted to $(-1,0)$ or $(0,1)$. Moreover,
\begin{eqnarray}\label{0368}  &&
\|\tilde{\varphi}\|_{L^2(-1,1)}^2+\|\tilde{\theta}\|_{L^2(-1,1)}^2+\|\tilde{\psi}\|_{L^2(-1,1)}^2=1,
\\
&& \label{0369} \|\tilde{\psi}'\|_{L^2(-1,1)}\leq {|\xi|\sqrt{2\varrho_+/{\varrho_-}}}.
\end{eqnarray}
\end{thm}
\begin{pf} In view of Theorem \ref{pro:0301}, we have a solution
$(\psi ,\lambda ):=(\psi (|\xi_1|,|\xi_2|,x_3),\lambda (|\xi_1|,|\xi_2|))$
satisfying (\ref{0309})--(\ref{0313}). Moreover, $\lambda >0$ and
$\psi\in\mathbb{A}$ is smooth when restricted to $(-1,0)$ or $(0,1)$.
Then, multiplying (\ref{0306})$_1$ and (\ref{0306})$_2$ by $\xi_1$
and $\xi_2$ respectively, adding the resulting equations, and
utilizing (\ref{0306})$_4$, we find that $\pi$ can be expressed by $\psi$, i.e.,
\begin{equation}\label{0370}\pi =\pi(|\xi_1|,|\xi_2|,x_3)
=[{\lambda\mu \psi'''-(\lambda^2\varrho+\lambda\mu
|\xi|^2+B^2\xi_1^2)\psi'}]/({\lambda|\xi|^{2}}).
\end{equation}

Notice that (\ref{0306})$_1$ can be rewritten as
\begin{equation}\label{0371}\varphi''-(\lambda^2\varrho+\lambda\mu |\xi|^2+B^2\xi_1^2)
\varphi/{(\lambda\mu)}=-{\xi_1\pi}/{\mu}  \end{equation}
with jump and boundary conditions
\begin{equation}\label{0372}
\llbracket  \varphi \rrbracket=0,\  \llbracket \mu
(\xi_1\psi-\varphi') \rrbracket=0,\ \varphi(-1)=\varphi(1)=0.
\end{equation}
Hence, we can easily construct a unique solution of the form
\begin{equation}\label{0373}\varphi =(\xi,x_3)=\left\{
                              \begin{array}{ll}
\xi_1(c_1e^{a_+x_3}+c_2e^{-a_+x_3}-f_+(x_3)),&\mbox{ on }(0,1),\\[0.4em]
\xi_1(c_3e^{a_-x_3}+c_4e^{-a_-x_3}-f_-(x_3)),&\mbox{ on }(-1,0)
\end{array}\right.
\end{equation}
 to the equation (\ref{0371}) with jump and boundary conditions (\ref{0372}),
where
\begin{eqnarray*} &&
a_\pm=\sqrt{|\xi|^2+\frac{\lambda \varrho}{\mu_\pm}+\frac{B^2\xi_1^2}{\lambda\mu_\pm}}, \\[1mm]
&&  f_\pm(x_3)=\frac{1}{2a_\pm\mu_\pm}\int_0^{x_3} \pi
(e^{a_\pm(x_3-y)}-e^{a_\pm(y-x_3)})\mathrm{d}y,
\end{eqnarray*}
and
\begin{equation*}
\left[
 \begin{array}{c}
c_1 \\
c_2 \\
c_3\\c_4
\end{array}
\right] =\left[
 \begin{array}{cccc}
1 & 1 & -1&-1 \\
\mu_+ a_+& -\mu_+a_+ & -\mu_-a_- & \mu_-a_- \\
e^{a_+}& e^{-a_+} & 0 &0 \\
0 & 0& e^{-a_-}& e^{a_-} \end{array}
 \right]^{-1}\left[
\begin{array}{c}
0\\
(\mu_+-\mu_-)\psi(0) \\
f(1) \\
f(-1)
\end{array}
\right].
\end{equation*}

Similarly to (\ref{0373}),
\begin{equation}\theta:=\theta(\xi,x_3)=\left\{
                              \begin{array}{ll}
\xi_2(c_1e^{a_+x_3}+c_2e^{-a_+x_3}-f_+(x_3)),&\mbox{ on
}(0,1),\\[0.4em]
\xi_2(c_3e^{a_-x_3}+c_4e^{-a_-x_3}-f_-(x_3)),&\mbox{ on }(-1,0)
\end{array}\right.
\end{equation}
is a unique solution of (\ref{0306})$_2$ with jump and boundary
conditions:
\begin{equation*}
                             \llbracket
\theta\rrbracket=0,\;\; \llbracket \mu (\xi_2\psi-\theta')
\rrbracket=0,\;\; \theta(-1)=\theta(1)=0.
                        \end{equation*}

Consequently, $({\varphi},{\theta},{\psi},{\pi})$ is a solution to
the system (\ref{0306})--(\ref{0308}).  Now, we define
\begin{equation*}\begin{aligned}\label{}
(\tilde{\varphi},\tilde{\theta},\tilde{\psi},\tilde{\pi}):=
&(\tilde{\varphi}(\xi,x_3),\tilde{\theta}(\xi,x_3),\tilde{\psi}(\xi,x_3),\tilde{\pi}(\xi,x_3))\\
 =&                            {({\varphi},{\theta},{\psi},{\pi})}/
                            (
                            \|\varphi\|_{L^2(-1,1)}^2+\|\theta\|_{L^2(-1,1)}^2+\|\psi\|_{L^2(-1,1)}^2).
\end{aligned}\end{equation*}
Thus, $(\tilde{\varphi},\tilde{\theta},\tilde{\psi},\tilde{\pi})$ is
still a solution to the system (\ref{0306})--(\ref{0308}), and
moreover,
$(\tilde{\varphi},\tilde{\theta},\tilde{\psi},\tilde{\pi})$
satisfies (\ref{0368}).

Finally, making use of (\ref{0306})$_4$ and (\ref{0332}), we conclude that
$$ \begin{aligned}
\frac{1}{\varrho_+|\xi|^2}= &
\frac{1}{2\varrho_+|\xi|^2}\int_{-1}^1\varrho(|\xi|^2|\psi|^2+|\psi'|^2)\mathrm{d}x_3\\
\leq
&\int_{-1}^1(|\varphi|^2+|\theta|^2+|\psi|^2)\mathrm{d}x_3\\
= &
\|\varphi\|_{L^2(-1,1)}^2+\|\theta\|_{L^2(-1,1)}^2+\|\psi\|_{L^2(-1,1)}^2
\end{aligned} $$
 and
$$  \int_{-1}^1|\psi'|^2\mathrm{d}x_3\leq {2}/{\varrho_-}. $$
The above two inequalities imply (\ref{0369}) immediately.  \hfill
$\Box$
\end{pf}
\begin{rem}\label{rem:0301}
For each $x_3$, it is easy to see that the solution
$(\tilde{\varphi}(\xi,\cdot),\tilde{\theta}(\xi,\cdot),\tilde{\psi}(\xi,\cdot),\tilde{\pi}(\xi,\cdot),
\lambda(|\xi|))$ constructed in Theorem \ref{thm:0302} possesses the following properties:
\begin{enumerate}[\quad \ (1)]
  \item $\lambda(|\xi_1|,|\xi_2|)$, $\tilde{\psi}(\xi,\cdot)$ and $\tilde{\pi}(\xi,\cdot)$
  are even on $\xi_1$ or $\xi_2$, when the other variable is fixed;
  \item $\tilde{\varphi}(\xi,\cdot)$ is odd on $\xi_1$, but even on $\xi_2$,
  when the other variable is fixed;
  \item $\tilde{\theta}(\xi,\cdot)$ is even on $\xi_1$, but odd on $\xi_2$,
  when the other variable is fixed.
  \end{enumerate}
\end{rem}

The next result provides an estimate for the $H^k$-norm of the
solutions $(\varphi,\theta,\psi,\pi)$ with $\xi$ varying in domain
$\mathbb{D}$, which will be useful in the next section when such
solutions are integrated in a Fourier synthesis. To emphasize the
dependence on $\xi$, we write these solutions as
$\big(\tilde{\varphi}(\xi)=\tilde{\varphi}(\xi,x_3),\tilde{\theta}(\xi)
=\tilde{\theta}(\xi,x_3),\tilde{\psi}(\xi)=\tilde{\psi}(\xi,x_3),\tilde{\pi}(\xi)
=\tilde{\pi}(\xi,x_3)\big)$.

\begin{lem}\label{lem:0302}
Let $\xi\in \mathbb{D}$, $\theta(\xi):=\tilde{\theta}(\xi)$,
$\psi(\xi):=\tilde{\psi}(\xi)$, $\pi(\xi):=\tilde{\pi}(\xi)$ and
$\lambda(|\xi_1|,|\xi_2|)$ be constructed in Theorem
\ref{thm:0302}, then for any $k\geq  0$ there exit positive
constants $a_k$, $b_k$ and $c_k$ depending on $d_1$, $d_2$,
$\lambda_0$, $|B|$, $\varrho_\pm$, $\mu_\pm$ and $g$, so that
\begin{eqnarray}
&&\label{0382}
\|\varphi(\xi)\|_{H^k(-1,0)}+\|\varphi(\xi)\|_{H^k(0,1)}+\|\theta(\xi)\|_{H^k(-1,0)}
+\|\theta(\xi)\|_{H^k(0,1)}\leq a_k,\\
&&\label{0383}
\|\psi(\xi)\|_{H^k(-1,0)}+\|\psi(\xi)\|_{H^k(0,1)}\leq
b_k,\\
&&\label{0384} \|\pi(\xi)\|_{H^k(-1,0)}+\|\pi(\xi)\|_{H^k(0,1)}\leq
c_k,
\end{eqnarray}
where $d_1$ and $d_2$ are the same constants in (\ref{0320}), and
$\lambda_0$ is defined by (\ref{0367}). Moreover,
\begin{equation}\label{0385}
\|\varphi\|_{L^2(-1,1)}^2+\|\theta\|_{L^2(-1,1)}^2+\|\psi\|_{L^2(-1,1)}^2=1.
\end{equation}
\end{lem}
\begin{pf} Throughout this proof, we denote by $\tilde{c}_1,\cdots,\tilde{c}_8$ generic
positive constants which may depend on $d_1$, $d_2$, $\lambda_0$,
$|B|$, $\varrho_\pm$, $\mu_\pm$ and $g$, but not on $|\xi|$.
Obviously, (\ref{0385}) follows from (\ref{0369}) immediately.

(i) First we write (\ref{0309}) as
\begin{equation}\label{0386}
\psi''''(\xi)=\big[{(\lambda^2\varrho+2\lambda\mu
|\xi|^2+B^2\xi_1^2)\psi''(\xi)-|\xi|^2(\lambda^2\varrho +\lambda\mu
|\xi|^2+B^2\xi_1^2)\psi(\xi)}\big]/(\lambda\mu).
\end{equation}

If we make use of (\ref{0348}), (\ref{0367}), and the fact $|\xi|\leq
d_2$, Lemma \ref{lem:0301} and the Cauchy-Schwarz inequality, we see that
there exists two constants $\tilde{c}_1$ and $\tilde{c}_2$, such that
\begin{equation}\label{0387}\begin{aligned}
\|\psi''''(\xi)\|_{L^2(I_\pm)}&\leq \tilde{c}_1\big(\|\psi(\xi)\|_{L^2(I_\pm)}
+\|\psi''(\xi)\|_{L^2(I_\pm)}\big)\\
& \leq (\tilde{c}_2+1)\big(\varepsilon^{-1/2}\|\psi(\xi)\|_{L^2(I\pm)}
+\sqrt{\varepsilon}\|\psi''''(\xi)\|_{L^2(I_\pm)}\big)
\;\;\mbox{ for any }\varepsilon\in (0,1),
\end{aligned}\end{equation}
respectively, where $I_+=(0,1)$ and $I_-=(-1,0)$. Choosing
$\sqrt{\varepsilon}=1/\{2(\tilde{c}_2+1)\}$ in (\ref{0387}) and
using (\ref{0385}), we arrive at
\begin{equation}\label{0388}
\|\psi''''(\xi)\|_{L^2(I_\pm)}\leq \tilde{c}_3 \quad\mbox{for some
constant  }\tilde{c}_3>0,
\end{equation}
whence,
\begin{equation}\label{0389}
\|\psi''''(\xi)\|_{L^2(-1,1)}\leq \tilde{c}_4:=2\tilde{c}_3.
\end{equation}

Thus from (\ref{0388}) and Lemma \ref{lem:0301} we get
\begin{equation}\label{0391}
\|\psi''(\xi)\|_{L^2(-1,1)}\leq \tilde{c}_5.
\end{equation}
and
\begin{equation}\label{0394}
\|\psi'''(\xi)\|_{L^2(-1,1)}\leq \tilde{c}_6.
\end{equation}

Summarizing the estimates (\ref{0385}), (\ref{0369}), and (\ref{0389})--(\ref{0394}),
we conclude that, for each
nonnegative integer $k\in [0, 4]$, there exists a constant
$\tilde{b}_k>0$ depending on $d_2$, $|B|$, $\lambda_0$,
$\varrho_\pm$, $\mu_\pm$ and $g$, such that
\begin{equation}\label{0395}
\|\psi^{(k)}(\xi)\|_{L^2(-1,1)}\leq \tilde{b}_k.
\end{equation}

Differentiating (\ref{0386}) with respect to $x_3$ and using
(\ref{0395}), we find, by induction on $k$, that (\ref{0395}) holds
for any $k\geq 0$. This gives (\ref{0383}).

(ii) Recalling the expression (\ref{0370}) of $\pi$ and the fact
that $|\xi|\geq d_1$, we employ (\ref{0367}) and (\ref{0348}) to
deduce that for any $k\geq 0$,
$$ \begin{aligned}
\|\pi^{(k)}(\xi)\|_{L^2{(-1,1)}}\leq &
\frac{\mu^+}{|\xi|^2}\|\psi^{(k+3)}(\xi)\|_{L^2{(-1,1)}}
+\left(\frac{\lambda\varrho^+}{|\xi|^2}+\mu^++\frac{B^2\xi_1^2}{\lambda|\xi|^2}
\right)\|\psi^{(k+1)}(\xi)\|_{L^2{(-1,1)}}  \\
\leq &\frac{\mu^+}{d_1^2}\tilde{b}_{k+3}
+\left(\frac{\Lambda\varrho^+}{d_1^2}+\mu^++\frac{B^2}{\lambda_0}\right)\tilde{b}_{k+1},
\end{aligned} $$
which implies (\ref{0384}).

(iii) Making use of (\ref{0371}), (\ref{0384}), (\ref{0385}),
(\ref{0348}) and $|\xi|\leq d_2$, we get
\begin{equation}\label{0397}\begin{aligned}
\|\varphi{''}(\xi)\|_{L^2(-1,1)}\leq&\frac{d_2}{\mu_-}\|\pi(\xi)\|_{L^2(-1,1)}+
\left(\frac{\Lambda\varrho^+}{\mu_-}+
d_2^2+\frac{B^2d_2^2}{\lambda_0\mu_-}\right)\|\varphi(\xi)\|_{L^2(-1,1)}
\leq \tilde{c}_{7}.
\end{aligned}\end{equation}
Applying (\ref{0397}), (\ref{0385}) and Lemma \ref{lem:0301}, we obtain
\begin{equation}\label{0398}\begin{aligned}
\|\varphi{'}(\xi)\|_{L^2(-1,0)}+\|\varphi{'}(\xi)\|_{L^2(0,1)}\leq&
\tilde{c}_{8}.
\end{aligned}\end{equation}
Combining (\ref{0385}) with (\ref{0357}) and (\ref{0398}), we
conclude that, for each nonnegative integer $k\in [0, 2]$, there
exists a constant $\tilde{a}_k>0$ depending on $d_1$, $d_2$,
$\lambda_0$, $\varrho_\pm$, $\mu_\pm$ and $g$, so that
\begin{equation}\label{0399}
\|\varphi^{(k)}(\xi)\|_{L^2(-1,1)}\leq \tilde{a}_k.
\end{equation}

Thus, by virtue of (\ref{0371}) and induction on $k$, (\ref{0399})
holds for any $k\geq 0$. Following the same procedure as used in
estimating $\varphi$, we infer that for each $k\geq 0$,
\begin{equation}\label{03100} \|\theta^{(k)}(\xi)\|_{L^2(-1,1)}\leq
\tilde{d}_k
\end{equation}
for some constant $\tilde{d}_k$  depending on $d_1$, $d_2$,
$\lambda_0$, $\varrho_\pm$, $\mu_\pm$ and $g$. Adding (\ref{03100})
to (\ref{0399}), we arrive at
\begin{equation*}\|\varphi^{(k)}(\xi)\|_{L^2(-1,1)} +\|\theta^{(k)}(\xi)\|_{L^2(-1,1)}\leq
(\tilde{a}_k+\tilde{d}_k) \mbox{ for any }k\geq 0,
\end{equation*}
which yields (\ref{0382}). This completes the proof. \hfill $\Box$
\end{pf}

\subsection{Fourier synthesis}

In this section we will use the Fourier synthesis to build growing
solutions to (\ref{0134}) with $\bar{B}=(B,0,0)$ out of the
solutions constructed in the previous section (Theorem
\ref{thm:0302}) for fixed spatial frequency $\xi\in \mathbb{A}$. The
constructed solutions will grow in-time in the piecewise Sobolev
space of order $k$, $\mathrm{\bar{H}}^k(\Omega)$, defined by
(\ref{0203}).
\begin{thm}\label{thm:0303}
 Let $\bar{B}=(B,0,0)$ and $f\in C_0^\infty[0,\infty)$ be a real-valued
function. For any $\xi\in \mathbb{R}^2$, we define
\begin{equation*}\label{} {w}(\xi,
x_3)=-i\tilde{\varphi}(\xi,x_3)e_1-i\tilde{\theta}(\xi,x_3)e_2+\tilde{\psi}(\xi,
x_3)e_3,
\end{equation*}
where
\begin{equation*}\label{03102} (\tilde{\varphi}, \tilde{\theta}, \tilde{\psi}, \tilde{\pi})(\xi, x_3)=
\left\{   \begin{array}{ll}
          \mbox{the solutions provided by Theorem \ref{thm:0302} } &\quad \\
     with\ \lambda(|\xi_1|,|\xi_2|)>0,&\mbox{if }\xi\in \mathbb{D}, \\[1mm]
zero,\quad  and\
\lambda(|\xi_1|,|\xi_2|)\equiv\lambda_0,&\mbox{if }
\xi\in\!\!\!\!\!/\ \mathbb{D}.
                  \end{array}
                \right.
\end{equation*}
Denote
\begin{eqnarray}
&&\label{03103}\eta(t,x)=\frac{1}{4\pi^2}\int_{\mathbb{R}^2}
f(|\xi|)w(\xi,x_3)e^{\lambda(|\xi_1|,|\xi_2|)t}e^{ix'\xi}\mathrm{d}\xi ,\\
&&\label{03104}v(t,x)=\frac{1}{4\pi^2}\int_{\mathbb{R}^2}
\lambda(|\xi_1|,|\xi_2|)f(|\xi|){w}(\xi,x_3)e^{\lambda(|\xi_1|,|\xi_2|)t}e^{ix'\xi}\mathrm{d}\xi ,\\
&&\label{03105}q(t,x)=\frac{1}{4\pi^2}\int_{\mathbb{R}^2}
\lambda(|\xi_1|,|\xi_2|)f(|\xi|)\tilde{\pi}(\xi,x_3)e^{\lambda(|\xi_1|,|\xi_2|)t}e^{ix'\xi}\mathrm{d}\xi .
\end{eqnarray}
Then, $(\eta,v,q)$ is a real-valued solution to the linearized problem (\ref{0134})
 along with jump and boundary conditions (\ref{0135}) and (\ref{0136}).
For every $k\in\mathbb{N}$, we have the estimate
\begin{equation}\label{03106}\begin{aligned}
\|\eta(0)\|_{\mathrm{\bar{H}}^k}+\|v(0)\|_{\mathrm{\bar{H}}^k}+\|q(0)\|_{\mathrm{\bar{H}}^k}
 \leq \tilde{c}_k\Big(\int_{\mathbb{R}^2} f^2(|\xi|)\mathrm{d}\xi\Big)^{1/2}<\infty,
\end{aligned}
\end{equation}
where $\tilde{c}_k>0$ is a constant depending on the parameters
$d_1$, $d_2$, $|B|$, $\lambda_0$,  $\varrho_\pm$, $\mu_\pm$ and $g$.
Moreover, for every $t>0$ we have $\eta(t),v(t),q(t)\in
\mathrm{\bar{H}}^k(\Omega)$, and
\begin{eqnarray}
&&\label{03107} e^{t\lambda_0}\|\eta(0)\|_{\mathrm{\bar{H}}^k}
\leq \|\eta(t)\|_{\mathrm{\bar{H}}^k}\leq e^{t\Lambda}\|\eta(0)\|_{\mathrm{\bar{H}}^k},\\
&&\label{03108}e^{t\lambda_0}\|v(0)\|_{\mathrm{\bar{H}}^k}\leq
\|v(t)\|_{\mathrm{\bar{H}}^k}\leq
e^{t\Lambda}\|v(0)\|_{\mathrm{\bar{H}}^k},\\
&&\label{03109}e^{t\lambda_0}\|q(0)\|_{\mathrm{\bar{H}}^k}\leq
\|q(t)\|_{\mathrm{\bar{H}}^k}\leq
e^{t\Lambda}\|q(0)\|_{\mathrm{\bar{H}}^k},
\end{eqnarray} where $\lambda_0$ and $\Lambda$ are defined in
(\ref{0367}) and (\ref{0348}), respectively.
\end{thm}
\begin{pf}
For each fixed $\xi\in \mathbb{R}^2$,
\begin{eqnarray*}
&& \eta (t,x')=f(|\xi|){w}(\xi,x_3)e^{\lambda(|\xi_1|,|\xi_2|)t}e^{ix'\xi},\\[1mm]
&& v(t,x)=\lambda(|\xi_1|,|\xi_2|)f(|\xi|){w}(\xi,x_3)e^{\lambda(|\xi_1|,|\xi_2|)t}e^{ix'\xi}\\[1mm]
&&
q(t,x)=\lambda(|\xi_1|,|\xi_2|)f(|\xi|){\pi}(\xi,x_3)e^{\lambda(|\xi_1|,|\xi_2|)t}e^{ix'\xi}
\end{eqnarray*}
give a solution to (\ref{0134}). Since $f\in C_0^\infty[0,\infty)$,
by the construction of $w$, Lemma \ref{lem:0302} implies that
\begin{equation*}\label{}
\sup_{\xi\in\mathrm{supp}(f)}\|\partial_{3}^k{w}(\xi,\cdot)\|_{L^\infty(\Omega)}<\infty\quad
\mbox{ for all }k\in \mathbb{N}.\end{equation*} Also,
$\lambda(\xi)\leq \Lambda$. These bounds show that the Fourier
synthesis of the solution given by (\ref{03103})--(\ref{03105}) is
also a solution of (\ref{0134}). Because $f$ is real-valued and
radial, $\mathbb{D}$ is a symmetrical domain (see Proposition
\ref{pro:0301}), we can easily verify, recalling Remark
\ref{rem:0301}, that the Fourier synthesis is real-valued.

The estimate (\ref{03106}) follows from Lemma \ref{lem:0302} with
arbitrary $k\geq 0$ and the construction of ($\eta, v, q$). Finally,
we can use (\ref{0367}), (\ref{0348}) and
(\ref{03103})--(\ref{03105}) to obtain the estimates
(\ref{03107})--(\ref{03109}). \hfill $\Box$
\end{pf}

\section{Global instability for the linearized problem}

\subsection{Uniqueness of the linearized equations}
In this section, we will show the uniqueness of solutions to the
linearized problem with lower regularity, which will be applied to the proof of
Theorem \ref{thm:0202} in Section 5. For this purpose, we need
a generalized formula of integrating by parts (or Gauss-Green's formula).
Let us first recall the boundary trace theorem (see
Theorem 5.36 in \cite[Chaperter 5]{ARAJJFF}).
\begin{lem}\label{lem:0401}
Let $U$ be a domain in $\mathbb{R}^n$ satisfying the uniform
$C^m$-regularity condition, and assume that there exists a simple
($m$, $p$)-extension operator $E$ for $U$. Also assume that $mp<n$
and $p\leq q\leq p^*=(n-1)p/(n-mp)$. Then, there exists a bounded
linear operator
\begin{equation*}\gamma^{U}:~W^{m,p}(U)\rightarrow L^q(\partial U),
\end{equation*}
such that
\begin{equation*}\gamma^{U}(u)=u\mbox{ on }\partial U
\end{equation*} for all $u\in W^{m,p}(U)\cap C(\bar{U})$.
\end{lem}
The function $\gamma^U(u)\in L^q(\partial U)$ is called the trace of
the function of $u\in W^{1,p}(U)$ on the boundary $\partial U$. By
the Stein extension theorem (see Theorem 5.24 in \cite[Chaperter 5]{ARAJJFF})
and the definition of the uniform $C^m$-regularity
condition (see Definition 4.10 in \cite[Chaperter 5]{ARAJJFF}), it
is easy to verify that $\Omega$, $\Omega_+$ and $\Omega_-$ have
different simple ($m$, $p$)-extension operators. Keeping these facts
in mind, we can start to show the following formula of integrating
by parts. For convenience in the subsequent analysis, we will use
the notations $\gamma_+(f):=\gamma^{\Omega_+}(f_+)$ and
$\gamma_-(f):=\gamma^{\Omega_+}(f_-)$.
\begin{lem}
For all $v\in H^1_0(\Omega)$ and $w\in \mathrm{\bar{H}}^1(\Omega)$,
we have
\begin{equation}\label{0403}\int_\Omega \partial_iwv\mathrm{d}x=-\int_\Omega
w\partial_iv\mathrm{d}x+\int_{\mathbb{R}^2}
(\gamma_+(w)-\gamma_-(w))\gamma_+(v)\alpha_i\mathrm{d}x
\end{equation} for $i=1,2,3$, where $\alpha_1=\alpha_2=0$ and
$\alpha_3=-1$.
\end{lem}
\begin{pf} Temporarily suppose $\bar{v}\in C^1_0({\Omega})$, $\bar{w}_+\in C^1(\bar{\Omega}_+)$ and $\bar{w}_-\in
C^1(\bar{\Omega}_-)$. By the Gauss-Green theorem, we have
\begin{equation}\label{0404}
\begin{aligned}\int_\Omega \partial_i\bar{w}\bar{v}\mathrm{d}x=&-\int_\Omega
\bar{w}\partial_i\bar{v}\mathrm{d}x+\int_{\mathbb{R}^2}
((\bar{w}_+-\bar{w}_-)\bar{v})(x',0)\alpha_i\mathrm{d}x.
\end{aligned}\end{equation}

Using Lemma \ref{lem:0401}, one has
$$ \begin{aligned}
\|(\bar{v}-\gamma_+(v))(x',0)\|_{L^2({\mathbb{R}^2})}\leq&
\|\bar{v}-\gamma_+(v)\|_{L^2(\partial \Omega_+)}
\\ =&\|\gamma_+(\bar{v}-v)\|_{L^2(\partial
\Omega_+)}\leq c
\|\bar{v}-v\|_{H^1(\Omega_+)}\\
\leq &c \|\bar{v}-v\|_{H^1_0(\Omega)}\end{aligned} $$
 and
$$ \| (\bar{w}_+ -\gamma_+(w_+))(x',0)\|_{L^2({\mathbb R}^2)} \leq
c\|\bar{w}_+ - w_+\|_{H^1(\Omega_+)} $$
for some constant $c>0$.  By the H\"{o}lder inequality, the above
two estimates imply that
\begin{equation}\label{0407}\begin{aligned}
& \|(\bar{w}_+\bar{v}-\gamma_+(w)\gamma_+(v))(x',0)\|_{L^1({\mathbb{R}^2})}\\
&
\leq\|(\bar{v}(\bar{w}_+-\gamma_+({w})))(x',0)\|_{L^1({\mathbb{R}^2})}+
\|(\gamma_+(w)(\bar{v}-\gamma_+(v)))(x',0)\|_{L^1({\mathbb{R}^2})}\\
&
\leq\|\bar{v}(x',0)\|_{L^2({\mathbb{R}^2})}\|(\bar{w}_+-\gamma_+({w}))(x',0)\|_{L^2({\mathbb{R}^2})}
\\ & \quad
+\|\gamma_+(w)(x',0)\|_{L^2({\mathbb{R}^2})}\|(\bar{v}-\gamma_+(v))(x',0)\|_{L^2({\mathbb{R}^2})}\\
& \leq c^2\|\bar{v} \|_{H^1_0(\Omega)}\|\bar{w}_+-w_+
\|_{H^1(\Omega_+)} +c^2\|w_+ \|_{H^1(\Omega_+)}\|\bar{v}-v
\|_{H^1_0(\Omega)}.\end{aligned}
\end{equation}

Similarly to (\ref{0407}), one gets
\begin{equation}\label{0408}\begin{aligned} &
\|(\bar{w}_-\bar{v}-\gamma_-({w})\gamma_+(v))(x',0)\|_{L^1({\mathbb{R}^2})}\\
& \quad \leq c^2(\|\bar{v}
\|_{H^1_0(\Omega)}\|\bar{w}_--w_-\|_{H^1(\Omega_-)}
+\|w_-\|_{H^1(\Omega_-)}\|\bar{v}-v\|_{H^1_0(\Omega)}).
\end{aligned}\end{equation}
In addition, if $\bar{v}_{m}\rightarrow v$ strongly in
$H_0^1(\Omega)$, then there exists $m_0>0$ such that
\begin{equation}\label{0409}\begin{aligned}
\|\bar{v}_m\|_{H^1_0(\Omega)}\leq \|v\|_{H^1_0(\Omega)}+1\qquad
\mbox{ for any }m\geq m_0.
\end{aligned}\end{equation}

Since $C_0(\Omega)$ is dense in $H_0^1(\Omega)$ and
$C_0(\mathbb{R}^3)$ dense in $H^{1}(\Omega_\pm )$, the identity (\ref{0403}) follows from
(\ref{0404})--(\ref{0409}) and a standard density argument.
\hfill $\Box$
\end{pf}
\begin{definition}\label{def:0401}
Given $T>0$ and the initial data ($\eta_0,v_0$) to the linearized
problem (\ref{0134})--(\ref{0136}), a triple ($\eta,v,q$) is called
a strong solution of (\ref{0134})--(\ref{0136}), if
\begin{enumerate}[\quad \ (1)]
  \item $\eta$, $v\in C^0([0,T],L^2(\Omega))$, $\eta(0)=\eta_0$, $v(0)=v_0$ and
 \begin{equation}\label{0410}\mathrm{ess}\sup_{0< t<T}
 (\|{v}(t)\|_{\mathrm{\bar{H}}^3}+\|{\eta}(t)\|_{\mathrm{\bar{H}}^3}
 +\|{q}(t)\|_{\mathrm{\bar{H}}^1}+\|{v}(t)\|_{H^1_0(\Omega)}) < \infty .
 \end{equation}
 \item The equations
  \begin{eqnarray}\label{0411} &&
\partial_t\eta=v, \\   \label{0412}
&&\varrho \partial_t v+\nabla q=\mu \Delta v+\sum_{1\leq l,m\leq 3}
\bar{B}_l\bar{B}_m\partial_{lm}^2\eta, \\
\label{0413}
&& \; \mathrm{div}\, v=0   \end{eqnarray} hold a.e. in $(0,T]\times
(\Omega\setminus\{x_3=0\})$.
\item For a.e. $t\in (0,T)$,
\begin{eqnarray}
&& \int_{\mathbb{R}^2}(\gamma_+(q)I-\mu_+(\nabla v_++\nabla
v^T_+)-(\gamma_-(q)I-\mu_-(\nabla v_-+\nabla
v_-^T))) e_3\cdot\varphi\mathrm{d}x' \nonumber  \\
&& =\int_{\mathbb{R}^2}
\left(g[\varrho]\eta_{+3}\varphi_3+\sum_{l=1}^3\bar{B}_3\bar{B}_l(
\partial_l \eta_+ -\partial_l \eta_-) \cdot\varphi-\kappa\sum_{i=1}^2
\partial_i\eta_{+3}\partial_i\varphi_3\right)\mathrm{d}x'
\label{0414}
\end{eqnarray}
holds for any $\varphi\in H_0^1(\mathbb{R}^2)$, where
$\eta_{+3}$ and ${v}_3$ are the third component of $\eta_+$ and $v$,
respectively.
\end{enumerate}
\end{definition}
\begin{rem}\label{rem:0401}
Since $v(t)\in H^1_0(\Omega)\cap \mathrm{\bar{H}}^3(\Omega)$ for
each $t\geq 0$, we can make use of the embedding theorem and
(\ref{0413}) to obtain
\begin{eqnarray} &&
v(t)\in  C^0(\bar{\Omega}),\; v_+(t)\in C^1(\bar{\Omega}_+),\; v_-(t)\in C^1(\bar{\Omega}_-),
\label{04161111} \\
&&\label{0416} v(t)\equiv 0\mbox{ on }\partial{\Omega}, \\
 && \label{0417}
\nabla_{x'}v_+\equiv\nabla_{x'}v_-\mbox{ on }\mathbb{R}^2,\\
&& \mathrm{div}\, v(t)\equiv 0\;\;\mbox{ in }\bar{\Omega}\quad\mbox{for a.e. } t\geq 0.
\nonumber
\end{eqnarray}
 Thus, in view of (\ref{0417}), we define for the sake of simplicity that
$$ \nabla_{x'} v:=\nabla_{x'}v_+=\nabla_{x'}v_- \;\;\mbox{ on }\mathbb{R}^2\times\{0\}. $$
Moreover, by virtue of Lemma \ref{lem:0401}, there is a constant $c$
such that
\begin{equation} \label{0420}
\|v(t,x',0)\|_{H^1(\mathbb{R}^2)}\leq c\|{v}
(t)\|_{H^2(\Omega_\pm)}\quad\mbox{ for  a.e. }t\geq
0.\end{equation}
\end{rem}
\begin{rem}\label{rem:0403} In view the regularity of $(\eta,v,q)$ in Definition \ref{def:0401},
we can see that the equality (\ref{0414}) makes sense.
\end{rem}
\begin{rem}\label{rem:0402} It is easy to verify that any $(\eta$, $v$, $q$), which is
a solution established in Theorem \ref{thm:0303}, is a strong
solution to the linearized system (\ref{0134})--(\ref{0136}).
\end{rem}
\begin{thm} {\rm (Uniqueness)} \label{thm:0401}
Let $\bar{B}:=(\bar{B}_1,\bar{B}_2,\bar{B}_3)$ be a constant vector.
Assume that $(\tilde{\eta}, \tilde{v}, \tilde{q})$ and $(\bar{\eta},
\bar{w}, \bar{q})$ are two strong solutions to
(\ref{0134})--(\ref{0136}), with $\tilde{v}(0)=\bar{v}(0)=v_0$,
$\tilde{\eta}(0)=\bar{\eta}(0)=\eta_0$. Then
$(\tilde{\eta},\tilde{v},\nabla
\tilde{q})=(\bar{\eta},\bar{v},\nabla \bar{q})$.
\end{thm}
\begin{pf}
Let $(\eta,v, q)=(\tilde{\eta}-\bar{\eta},\tilde{v}-\bar{v},
\tilde{q}-\bar{q})$. Recalling Definition \ref{def:0401}, $(\eta,v,
q)$ is still a strong solution to the linearized system
(\ref{0134})--(\ref{0136}) with zero initial data, i.e., $\eta(0)=0$ and $v(0)=0$.

Let $t_1\in (0,T]$. Multiplying (\ref{0412}) by $v$, integrating
over $(0,\tau)\times\Omega$ for any $\tau\in (0,t_1)$ and using
(\ref{0412}), we find that
\begin{equation}\label{0421}\begin{aligned}
&\int_0^{\tau}\int_\Omega \varrho\partial_tv\cdot v\mathrm{d}x\mathrm{d}t
+\int_0^\tau\int_\Omega\mathrm{div}(qI-\mu(\nabla v +\nabla
v^T))\cdot v\mathrm{d}x\mathrm{d}t  \\
& \qquad =\sum_{1\leq l, m\leq 3}
\int_0^\tau\int_\Omega \bar{B}_l\bar{B}_m\partial_{lm}^2\eta\cdot
v\mathrm{d}x\mathrm{d}t.
\end{aligned}\end{equation}

(1) Firstly, we transform (\ref{0421}) to the form of
energy equality. By virtue of the regularity (\ref{0410}), (\ref{0412})
implies that
$$
\partial_t v\in L^2((0,T)\times \Omega), $$
which, together with $v\in L^\infty (0,T;H^1(\Omega))\cap C^0([0,T],L^2(\Omega))$,
yields
\begin{equation}\label{0422}\begin{aligned}
\int_0^{\tau}\int_\Omega \varrho
\partial_tv\cdot v\mathrm{d}x\mathrm{d}t =\frac{1}{2}\int_\Omega \varrho
v^2(\tau)\mathrm{d}x-\frac{1}{2}\int_\Omega \varrho
v^2(0)\mathrm{d}x.
\end{aligned}\end{equation}

Thanks to Lemma \ref{lem:0401}, (\ref{04161111})--(\ref{0416}), and
regularity of $q$, we obtain
\begin{equation}\label{0423}\begin{aligned}
&\int_0^\tau\int_\Omega\mathrm{div}(qI-\mu(\nabla v +\nabla v^T))
\cdot v\mathrm{d}x\mathrm{d}t\\
& =\int_0^\tau\int_{\mathbb{R}^2}(\gamma_-(q)I-\mu_-(\nabla
v_-+\nabla v_-^T)-(\gamma_+(q)I-\mu_+(\nabla
v_++\nabla v^T_+))) e_3\cdot v\mathrm{d}x'\mathrm{d}t\\
& \quad +\int_0^\tau\int_\Omega \mu\nabla v : (\nabla v+\nabla
v^T)\mathrm{d}x\mathrm{d}t.
\end{aligned}\end{equation}

Employing (\ref{0420}), $\mathrm{div}\,v =0$ at the plan $\{x_3=0\}$
 and the arguments similar to those used for (\ref{0403}), we can show
 \begin{equation*}\begin{aligned}
\sum_{j=1}^3\int_0^\tau\int_{\mathbb{R}^2}\mu v_j\partial_j v_{\pm3}
\mathrm{d}x'\mathrm{d}t=0.
\end{aligned}\end{equation*}
Using the above equality,  Lemma \ref{lem:0401}, and the
regularity of $u$ stated in Remark \ref{rem:0401}, we conclude that
 \begin{equation}\label{0424}\begin{aligned}
\int_0^\tau\int_\Omega\mu\nabla v:\nabla
v^T\mathrm{d}x\mathrm{d}t=\sum_{1\leq i, j\leq
3}\int_0^\tau\int_\Omega\mu\partial_i
v_j\partial_jv_i\mathrm{d}x\mathrm{d}t\equiv0.
\end{aligned}\end{equation}

In view of (\ref{0421})--(\ref{0424}) and $v(0)=0$, we find the
following energy equality
\begin{equation} \label{0425}\begin{aligned}
&\frac{1}{2}\int_\Omega \varrho v^2(\tau)\mathrm{d}x
+\int_0^\tau\int_\Omega \mu\nabla v : \nabla v
\mathrm{d}x\mathrm{d}t-\sum_{1\leq l, m\leq 3}\int_0^\tau\int_\Omega
\bar{B}_l\bar{B}_m\partial_{lm}^2\eta\cdot v\mathrm{d}x\mathrm{d}t \\
 &=\int_{\mathbb{R}^2}
(g[\varrho]\eta_{+3}v_3+\bar{B}_3\bar{B}_l(
\partial_l \eta_+-\partial_l \eta_-)\cdot v-\kappa\sum_{ i=1}^2\partial_i\eta_{+3}\partial_iv_3)\mathrm{d}x'.
\end{aligned}\end{equation}

Now we continue to transform the above energy equality by replacing
$\eta$ with $v$. Since $\eta\in C^0([0,T],L^2(\Omega))$ and
$\eta(0)=0$, the equation (\ref{0411}) gives
\begin{equation}\label{0426}
\eta(t,x)=\int_0^tv(s,x)\mathrm{d}s\mbox{ for any }t\geq 0,
\end{equation}
which, combined with (\ref{0417}), yields
\begin{equation}\label{0427}\begin{aligned}
\partial_i\eta(t,x):=\partial_i\eta_+(t,x)=\partial_i\eta_-(t,x)
=\int_0^t\partial_iv(s,x)\mathrm{d}s,\;\;\; i=1\mbox{ or }2
.\end{aligned}\end{equation} Using (\ref{0426}), (\ref{0427}), Lemma
\ref{lem:0401} and the regularity of $(\eta,v)$, we deduce that
\begin{equation}\label{0428}\begin{aligned}
\sum_{1\leq l,m\leq 3}\int_0^\tau\int_\Omega
\partial_{lm}^2\eta\cdot
v\mathrm{d}x\mathrm{d}t=&-\sum_{1\leq l,m\leq 3}\int_0^\tau\int_{\Omega}
\partial_m\eta\cdot\partial_l v\mathrm{d}x\mathrm{d}t-\int_0^\tau\int_{\mathbb{R}^2}
\llbracket\partial_3 \eta \rrbracket\cdot v\mathrm{d}x'\mathrm{d}t
\\
=&-\sum_{1\leq l,m\leq 3}\int_0^\tau\int_{\Omega}
\int_0^t\partial_mv(s,x)\mathrm{d}s\cdot\partial_l
v(t,x)\mathrm{d}x\mathrm{d}t\\
\qquad &-\int_0^\tau\int_{\mathbb{R}^2} \llbracket\partial_3 \eta
\rrbracket\cdot v\mathrm{d}x'\mathrm{d}t,
\end{aligned}\end{equation}
and
\begin{equation}\label{0429}
\sum_{i=1}^2\int_0^\tau\int_{\mathbb{R}^2}\partial_i\eta_{+3} \partial_i
v_3\mathrm{d}x'\mathrm{d}t =\sum_{i=1}^2\int_0^\tau\int_{\mathbb{R}^2}
\int_0^t\partial_i v_3(s,x',0)\mathrm{d}s\partial_i v_3(t,x',0)\mathrm{d}x'\mathrm{d}t.
\end{equation}

Consequently, inserting (\ref{0427})--(\ref{0429}) into (\ref{0425}), we arrive at
\begin{equation}\label{0430}\begin{aligned}
& \frac{1}{2}\int_\Omega \varrho v^2(\tau)\mathrm{d}x
+\int_0^\tau\int_\Omega\mu\nabla v:\nabla v
\mathrm{d}x\mathrm{d}t\\
& \qquad+\sum_{1\leq l,m\leq
3}\bar{B}_l\bar{B}_m\int_0^\tau\int_{\Omega}
\int_0^t\partial_mv(s,x)\mathrm{d}s\cdot\partial_l
v(t,x)\mathrm{d}x\mathrm{d}t \\
& \quad = g[\varrho]\int_0^\tau\int_{\mathbb{R}^2}\int_0^t
v_3(s,x',0)\mathrm{d}sv_3(t,x',0)\mathrm{d}x'\mathrm{d}t\\
&\qquad-\kappa\sum_{i=1}^2\int_0^\tau\int_{\mathbb{R}^2}
\int_0^t\partial_iv_3(s,x',0)\mathrm{d}s
\partial_iv_3(t,x',0)\mathrm{d}x'\mathrm{d}t
\end{aligned}\end{equation}

(2) Secondly, we further simplify the energy equality to an inequality.
  With the help of the regularity of $\partial_iv_3$, the
property of absolutely continuous functions and the Fubini theorem,
we conclude that
\begin{equation}\label{0431}\begin{aligned}
\int_0^\tau\int_{\mathbb{R}^2}\int_0^t
\partial_iv_3(s,x',0)\mathrm{d}s
\partial_iv_3(t,x',0)\mathrm{d}x'\mathrm{d}t
=&\int_{\mathbb{R}^2}
\int_0^\tau\int_0^t\partial_iv_3(s,x',0)\mathrm{d}s
\partial_iv_3(t,x',0)\mathrm{d}t\mathrm{d}x'\\
=&\int_{\mathbb{R}^2}\int_0^\tau\frac{d}{dt}\left[\int_0^t
\partial_iv_3(s,x',0)\mathrm{d}s\right]^2
\mathrm{d}t\mathrm{d}x'\\
=&\int_{\mathbb{R}^2}\left[\int_0^\tau
\partial_iv_3(t,x',0)\mathrm{d}t \right]^2\mathrm{d}x'\geq 0.
\end{aligned}\end{equation}
 Hence, by (\ref{0430})--(\ref{0431}), we find that
\begin{equation}\label{0433}\begin{aligned}
&\int_\Omega \varrho v^2(\tau)\mathrm{d}x+2\sum_{1\leq i,j\leq
3}^3\int_0^\tau\int_\Omega\mu|\partial_jv_i|^2\mathrm{d}x\mathrm{d}t
\\
& \quad \leq2g[\varrho]\int_0^\tau\int_{\mathbb{R}^2}\int_0^t
v_3(s,x',0)\mathrm{d}sv_3(t,x',0)\mathrm{d}x'\mathrm{d}t\\
&\qquad  -2\sum_{1\leq l, m\leq
3}\bar{B}_l\bar{B}_m\int_0^\tau\int_{\Omega}
\int_0^t\partial_mv(s,x)\mathrm{d}s\cdot\partial_l
v(t,x)\mathrm{d}x'\mathrm{d}t .
\end{aligned} \end{equation}

(3) Thirdly, we start to deduce the local-in-time uniqueness from the inequality
(\ref{0433}). Analogously to (\ref{0431}), the first
integral on the right-hand side of (\ref{0433}) can be bounded as follows.
\begin{equation}\label{0434}\begin{aligned}
& 2\int_0^\tau\int_{\mathbb{R}^2}\int_0^t
v_3(s,x',0)v_3(t,x',0)\mathrm{d}x'\mathrm{d}s\mathrm{d}t
=  \int_{\mathbb{R}^2}\left( \int_0^\tau v_3(t,x',0)\mathrm{d}t\right)^2\mathrm{d}x' \\
& \leq \tau
\int_0^\tau\int_{\mathbb{R}^2}|v_3(t,x',0)|^2\mathrm{d}x'\mathrm{d}t \\
& =2\tau \int_0^\tau\int_{\mathbb{R}^2}\int_1^0u_3(t,x',x_3)
\partial_3v_3(t,x',x_3)\mathrm{d}x_3\mathrm{d}x'\mathrm{d}t\\
& \leq \tau\int_0^\tau\int_{\mathbb{R}^2}\left (\frac{\mu}{\tau
g[\varrho]} \int_0^1|\partial_3v_3(t)|^2\mathrm{d}x_3+\frac{\tau
g[\varrho]}{\mu}
\int_0^1|v_3(t)|^2\mathrm{d}x_3\right)\mathrm{d}x'\mathrm{d}t\\
& \leq \frac{1}{g[\varrho]}\int_0^\tau
\|\sqrt{\mu}\partial_3v_3(t)\|_{L^2(\Omega)}^2\mathrm{d}t+\frac{\tau^2g[\varrho]}{\mu_-}
\int_0^\tau\|v(t)\|_{L^2(\Omega)}^2\mathrm{d}t.
\end{aligned}\end{equation}
On the other hand, letting
$B_0=\max\{|\bar{B}_1|,|\bar{B}_2|,|\bar{B}_3|,1\}$, and using the
H\"{o}lder, Minkowski and Cauchy-Schwarz inequalities, the second integral on the right-hand side of
(\ref{0433}) can be estimates as follows.
\begin{equation}\label{0435}\begin{aligned}
&\sum_{1\leq l, m\leq 3}\bar{B}_l\bar{B}_m\int_0^\tau\int_{\Omega}
\int_0^t\partial_mv(s,x)\mathrm{d}s\cdot\partial_l
v(t,x)\mathrm{d}x\mathrm{d}t\\
&\quad  \leq B_0^2\sum_{1\leq l,m\leq
3}\sum_{i=1}^3\int_0^\tau\left[\int_{\Omega}
\left(\int_0^t\partial_mv_i(s,x)\mathrm{d}s\right)^2\mathrm{d}x\right]^{1/2}
\left(\int_{\Omega}|\partial_l
v_i(t,x)|^2\mathrm{d}x\right)^{1/2}\mathrm{d}t\\
&\quad  \leq B_0^2\sum_{1\leq l, m\leq 3}
\sum_{i=1}^3\int_0^\tau\left[\int_0^t
\left(\int_{\Omega}|\partial_mv_i(s,x)|^2\mathrm{d}x\right)^{1/2}\mathrm{d}s\right]
\left(\int_{\Omega}|\partial_l
v_i(t,x)|^2\mathrm{d}x\right)^{1/2}\mathrm{d}t\\&\quad \leq
B_0^2\sum_{1\leq l,m\leq 3} \sum_{i=1}^3\left[\int_0^\tau
\left(\int_{\Omega}|\partial_mv_i(s,x)|^2\mathrm{d}x\right)^{1/2}\mathrm{d}s\right]\left[\int_0^\tau
\left(\int_{\Omega}|\partial_l
v_i(t,x)|^2\mathrm{d}x\right)^{1/2}\mathrm{d}t\right]\\
&\quad  \leq {{\tau}}B_0^2\sum_{1\leq i, j\leq 3}\int_0^\tau
\int_{\Omega}|\partial_jv_i(t,x)|^2\mathrm{d}x\mathrm{d}t.
\end{aligned}\end{equation}

Substituting (\ref{0434})--(\ref{0435}) into (\ref{0433}), we deduce that
\begin{equation}\label{0436}\begin{aligned}
&\|\sqrt{\varrho}v(\tau)\|^2_{L^2(\Omega)} +{\mu_-}\sum_{1\leq
i,j\leq
3}\int_0^\tau\|\partial_jv_i\|^2_{L^2(\Omega)}\mathrm{d}t\\
&\quad\leq {\tau^2g^2[\varrho]^2}{\mu_-^{-1}}
\int_0^\tau\|v(t)\|_{L^2(\Omega)}^2\mathrm{d}t
+{{\tau}}B_0^2\sum_{1\leq i, j\leq 3}\int_0^\tau
\|\partial_jv_i(t,x)\|^2_{L^2(\Omega)}\mathrm{d}t.
\end{aligned}
\end{equation}
Now, taking $t_1= \min\Big\{\frac{\mu_1}{2B_0^2},T\Big\}$,
the inequality (\ref{0436}) implies
\begin{equation}\label{0438}\|v(\tau)\|^2_{L^2(\Omega)}\leq
\frac{t_1^2g^2[\varrho]^2}{\mu_-\varrho_-}
\int_0^\tau\|v(t)\|_{L^2(\Omega)}^2\mathrm{d}t,\qquad \tau\in (0,t_1].
\end{equation}
Moreover, if we apply the Grownwall inequality to (\ref{0438}), we see that
\begin{equation*}
\|v(\tau)\|^2_{L^2(\Omega)}=0\quad\mbox{ for any }\tau\in [0,t_1],
\end{equation*}
which yields $v=0$, i.e., $\tilde{v}=\bar{v}$. This, combined with
(\ref{0412}) and (\ref{0426}), proves that
\begin{equation}\label{0437}(\tilde{\eta},\tilde{v},\nabla
\tilde{q})=(\bar{\eta},\bar{v},\nabla \bar{q})\quad\mbox{ for any
}\tau\in (0,t_1].\end{equation}

(4) Finally, the local-in-time uniqueness of the solution $(\tilde{\eta},\tilde{v})$
can be continued onto $[0,T]$. In fact, if $t_1<T$, we can take
$(\tilde{\eta},\tilde{v})(t_1)=(\bar{\eta},\bar{v})(t_1)$ as the
initial data, and continue the above procedure (1)--(3) to obtain
the uniqueness of solutions $(\tilde{\eta},\tilde{v})$ on $(0,t_2]$,
where $t_2=2t_1={\mu_1^2}/{(2b_0^2)}$ or $t_2=T$. Hence, after repeating
this procedure of extending time interval by finite times,
 (\ref{0437}) holds for any $\tau\in (0,T]$.  \hfill $\Box$
\end{pf}

\subsection{Proof of Theorem \ref{thm:0201}}
We define
\begin{equation*}
\beta_1=d_1+(d_2-d_1)/3,\quad \beta_2=d_1+ 2(d_2-d_1)/3.
\end{equation*}
Fix $j\geq k\geq 0$, $\alpha>0$ and let $\tilde{c}_j$ be the constants
from Theorem \ref{thm:0303}. For each $n\in \mathbb{N}$, let $t_n$ satisfy
\begin{equation}\label{0440}
{{e}^{2t_n\lambda_0}}= \alpha^2n^2\tilde{c}_j^2, \end{equation}
 i.e.,
\begin{equation*}\mbox{
}t_n=\frac{\mathrm{ln}\tilde{c}_j}{\lambda_0}+
\frac{\mathrm{ln}(\alpha
n)}{\lambda_0}:=C_{jk}+C_1\mathrm{ln}(\alpha n),
\end{equation*}
where $\lambda_0$ is defined by (\ref{0367}).
Choose $f_n\in C_0^\infty(\mathbb{R}^2)$, such that
supp$(f_n)\subset B(0, \beta_2)\backslash B(0,\beta_1)$, where $f_n$ is
real-valued and radial, and
\begin{equation}\label{0441}
\int_{\mathbb{R}^2}f_n^2(|\xi|)\mathrm{d}\xi=\frac{1}{\tilde{c}_j^2n^2}.
\end{equation}

Now, we can apply Theorem \ref{thm:0303} with $f=f_n$ to find that
$\big(\eta_n(t),v_n(t)$, $q_n(t)\big)\in \mathrm{\bar{H}}^j(\Omega)$
solves the problem (\ref{0134})--(\ref{0136}). It follows thus from
(\ref{03106}) and (\ref{0441}) that (\ref{0205}) holds for all $n$.

Recalling supp$(f_n)\subset B(\beta_1, \beta_2)$ and
$\lambda(|\xi_1|,|\xi_2|)\geq \lambda_0$, we have, after a
straightforward calculation and using (\ref{03104}), (\ref{0441}),
(\ref{0440}) and (\ref{0385}), that
\begin{equation*}
\begin{aligned}
\|v_n(t)\|_{\mathrm{\bar{H}}^k}^2  \geq&
\int_{\mathbb{R}^2}(1+|\xi|^2)^ke^{2t\lambda(|\xi_1|,|\xi_2|)}
f_n^2(|\xi|)\|{w}(\xi,x_3)\|_{L^2{(-1,1)}}^2\mathrm{d}\xi   \\
\geq& e^{2t\lambda_0}\int_{\mathbb{R}^2}f_n^2(|\xi|)\|{w}(\xi,x_3)\|_{L^2{(-1,1)}}^2\mathrm{d}\xi  \\
=&  e^{2(t-t_n)\lambda_0}\alpha^2
n^2\tilde{c}_j^2\int_{\mathbb{R}^2}
f_n^2(|\xi|)\mathrm{d}\xi  \\
\geq & \alpha^2\qquad \mbox{ for any }t\geq t_n,
\end{aligned}
\end{equation*}
which, together with $\eta_n(t,x)=\lambda(|\xi_1|,|\xi_2|)v_n(t,x)$,
implies (\ref{0206}) and (\ref{0207}). This completes the proof of
Theorem \ref{thm:0201}

\section{Proof of Theorem \ref{thm:0202}}
In this section we show Theorem \ref{thm:0202}. The main
idea of our proof comes from \cite{GYTI1,JFJSWWWO} but
with more complicated computations. We argue by contradiction. Therefore,
we suppose that the perturbed problem has the global stability of
order $k$ for some $k\geq 3$.

Let $\delta$, $C_2>0$ and $F:[0,\delta]\rightarrow \mathbb{R}^+$ be
the constants and function provided by the global stability of order
$k$, respectively. Fixing $n\in \mathbb{N}$ such that $n>C_2$.
Applying Theorem \ref{thm:0201} with this $n$, $t_n=T/2$, $k\geq 3$,
and $\alpha =2$, we find that
$(\tilde{\eta},\tilde{v},\tilde{\sigma})$ solves (\ref{0134}), satisfying
\begin{equation*}\label{0501}
\|\tilde{\eta}(0)\|_{\bar{\mathrm{H}}^k}
+\|\tilde{v}(0)\|_{\bar{\mathrm{H}}^k} <{n}^{-1}
\end{equation*}
 but
\begin{equation}\label{0502}
\|\tilde{v}(t)\|_{\bar{\mathrm{H}}^3}\geq 2\quad\mbox{ for }\; t\geq
T/2.
\end{equation}

For $\varepsilon>0$ we define
$\bar{\eta}_0^\varepsilon=\varepsilon\tilde{\eta}(0)$
 and $\bar{v}_0^\varepsilon=\varepsilon\tilde{v}(0)$.
Then, for $\varepsilon<\delta n$,
$\|(\bar{\eta}_0^\varepsilon,\bar{v}_0^\varepsilon)\|_{\bar{\mathrm{H}}^k}<\delta$.
So, according to the global stability of order $k$, there exist
${\eta}^\varepsilon,{v}^\varepsilon,{q}^\varepsilon$ that solve the
perturbed problem (\ref{0209})--(\ref{0212}) with $\bar{B}=(B,0,0)$, i.e.,
\begin{equation}\label{0503} \left\{
                              \begin{array}{ll}
\partial_t {\eta}^\varepsilon={v}^\varepsilon,\\
\varrho\partial_t
{v}^\varepsilon_i+(I_{jk}-{G^\varepsilon_{jk}})\partial_k{T}^\varepsilon_{ij}=B^2\partial_{11}^2
{\eta}^\varepsilon_i,\quad i=1, 2, 3,\\
\mathrm{div}{v}^\varepsilon- \mathrm{tr}(G^\varepsilon D
{v}^\varepsilon)=0
\end{array}
                            \right.
\end{equation}
 with jump conditions across the interface
\begin{equation}\label{0504}
\llbracket  {v}^\varepsilon \rrbracket =0,\ \llbracket
{T}_{ij}^\varepsilon {n}_j^\varepsilon
\rrbracket=g[\varrho]{\eta}_3^\varepsilon {n}_i^\varepsilon+{B}^2
(e_1+\partial_1\eta^\varepsilon)\cdot n^\varepsilon\llbracket
\partial_1\eta^\varepsilon_i \rrbracket+\kappa H^\varepsilon
n_i^\varepsilon\end{equation}
and initial data satisfying
$\|(\bar{\eta}_0^\varepsilon,\bar{v}_0^\varepsilon)\|_{\mathrm{\bar{H}}^k}<\delta$,
where
\begin{equation*}\label{0505}
(G^\varepsilon)^T:=(G_{jk}^\varepsilon)_{3\times3}^T:=I-(I+D\eta^\varepsilon)^{-1},\qquad
\qquad \qquad \qquad \qquad \qquad \qquad \
\end{equation*}
\begin{equation*}T_{ij}^\varepsilon=q^\varepsilon I_{ij}
-\mu\big( (I_{jk}-G_{jk}^\varepsilon)\partial_k
v_i^\varepsilon+(I_{ik}-G_{ik}^\varepsilon)\partial_k
v_j^\varepsilon)\big),\qquad \qquad \qquad \qquad\  \end{equation*}
\begin{equation}\label{0506}\begin{aligned}
n^\varepsilon=N^\varepsilon/|N^\varepsilon|
\mbox{ with }N^\varepsilon &= (e_1+\varepsilon\partial_1\bar{\eta}^\varepsilon)
\times (e_2+\varepsilon\partial_2\bar{\eta}^\varepsilon)\\
& = e_3+\varepsilon(e_1\times \partial_2
\bar{\eta}^\varepsilon+\partial_1\bar{\eta}^\varepsilon\times
e_2)+\varepsilon^2(\partial_1\bar{\eta}^\varepsilon\times\partial_2\bar{\eta}^\varepsilon)\\
 & =: e_3+\varepsilon\bar{N}^\varepsilon,
\end{aligned}\end{equation}
 and
\begin{equation}\label{0507}H^\varepsilon=\left(\frac{|e_1+\partial_1\eta^\varepsilon|^2\partial_2^2\eta^\varepsilon
-2(e_1+\partial_1\eta^\varepsilon)\cdot(e_2+\partial_2\eta^\varepsilon)\partial_1\partial_2\eta^\varepsilon+
|e_2+\partial_2\eta^\varepsilon|^2\partial_1^2
\eta^\varepsilon}{|e_1+\partial_1\eta^\varepsilon|^2|e_2+\partial_2\eta^\varepsilon|^2-|(e_1+\partial_1\eta^\varepsilon)
\cdot(e_2+\partial_2\eta^\varepsilon)|^2}\right)\cdot n^\varepsilon.
\end{equation}
Here we have used the Einstein convention of summing
over repeated indices, and $\eta^\varepsilon$ to denote the
both case $\eta_-^\varepsilon$ and $\eta_+^\varepsilon$ at the
interface $\{x_3=0\}$ in (\ref{0504}), (\ref{0506}) and
(\ref{0507}), except for the notation $\llbracket
\partial_l {\eta}_i^\varepsilon \rrbracket$ (see Remark \ref{rem:0201}).
 Moreover, ${\eta}^\varepsilon,{v}^\varepsilon,{q}^\varepsilon$ satisfy
\begin{equation}\label{0508}\sup_{0\leq t<T}\|({\eta}^\varepsilon,
{v}^\varepsilon,{q}^\varepsilon)(t)\|_{\bar{\mathrm{H}}^3} \leq
F(\|(\bar{\eta}_0^\varepsilon,\bar{v}_0^\varepsilon)\|_{\bar{\mathrm{H}}^3})
\leq C_2\varepsilon\|(\tilde{\eta}^\varepsilon,\tilde{v}^\varepsilon
)(0)\|_{\bar{\mathrm{H}}^3}<\varepsilon,\end{equation}
and $\eta_\pm^\varepsilon\in C^2(\bar{\Omega}_\pm)$ when $\kappa>0$.

Now, defining the rescaled functions
$\bar{\eta}^\varepsilon={\eta}^\varepsilon/\varepsilon$,
$\bar{v}^\varepsilon={v}^\varepsilon/\varepsilon$,
$\bar{q}^\varepsilon={q}^\varepsilon/\varepsilon$, and rescaling
(\ref{0508}), one gets
\begin{equation}\label{0509}
\sup_{0\leq t<T}\|(\bar{\eta}^\varepsilon,\bar{v}^\varepsilon,
\bar{q}^\varepsilon)(t)\|_{\bar{\mathrm{H}}^3} \leq 1,\end{equation}
which implies that
\begin{equation}\label{050911}
\sup_{0\leq t<T}\|(\bar{\eta}^\varepsilon,\bar{v}^\varepsilon,
\bar{q}^\varepsilon)(t)\|_{C^1(\bar{\Omega}_+)}+\sup_{0\leq
t<T}\|(\bar{\eta}^\varepsilon,\bar{v}^\varepsilon,
\bar{q}^\varepsilon)(t)\|_{C^1(\bar{\Omega}_-)} \leq K_1,
\end{equation}
for some constant $K_1$.  On the other hand, using the
jump conditions (\ref{0504}) in term of $v^\varepsilon$, one has
\begin{equation}\label{0510}
\sup_{0\leq t<T}\|\bar{v}^\varepsilon(t)\|_{{H}^1_0(\Omega)} \leq 1.
\end{equation}

Note that by construction, $(\bar{\eta}^\varepsilon,
\bar{v}^\varepsilon)(0)=(\tilde{\eta}, \tilde{v})(0)$. Our next goal
is to show that the rescaled functions converge as $\varepsilon\to 0$
to the solution
$(\bar{\eta},\bar{v},\bar{\sigma})$ of the linearized equations
(\ref{0134}) with $\bar{B}=(B,0,0)$ and the corresponding linearized
jump conditions.

\subsection{Convergence to the linearized equations}

We may further assume that $\varepsilon$ is sufficiently small, so that
\begin{equation}\label{0511}\sup_{0\leq t<T}\|\varepsilon D
\bar{\eta}^\varepsilon(t)\|_{L^\infty(\Omega)}<1/9 \;\;\mbox{ and }\;\;\varepsilon<1/(2K_2),
\end{equation}
where $K_2>0$ is the best constant in the inequality
\begin{equation*}\label{0512}\|FG\|_{\bar{\mathrm{H}}^2}\leq K_2
\|F\|_{\bar{\mathrm{H}}^2}\|G\|_{\bar{\mathrm{H}}^2}\end{equation*}
for $3\times 3$ matrix-valued functions $F$, $G$. Then,
\begin{equation*}\label{0513}
\bar{G}^\varepsilon=\big( I-(I+\varepsilon
(D\bar{\eta}^\varepsilon)^T)^{-1}\big)/\varepsilon\end{equation*} is
well-defined by (\ref{0511}) and the uniform boundedness in
$L^\infty (0,\infty; \bar{\mathrm{H}}^2(\Omega))$, since
\begin{equation}\label{0514}\begin{aligned}\|\bar{G}^\varepsilon\|_{\bar{\mathrm{H}}^2}
= &\left\|\sum_{n=1}^\infty(-\varepsilon
)^{n-1}(D\bar{\eta}^\varepsilon)^n\right\|_{\bar{\mathrm{H}}^2}\leq
\sum_{n=1}^\infty
\varepsilon^{n-1}\|(D\bar{\eta}^\varepsilon)^n\|_{\bar{\mathrm{H}}^2}\\
\leq &\sum_{n=1}^\infty (\varepsilon
K_1)^{n-1}\|D\bar{\eta}^\varepsilon\|^n_{\bar{\mathrm{H}}^2}\leq
\sum_{n=1}^\infty\frac{1}{2^{n-1}}\|\bar{\eta}^\varepsilon\|^n_{\bar{\mathrm{H}}^3}<\sum_{n=1}^\infty\frac{1}{2^{n-1}}=2.
\end{aligned}\end{equation}

Next, we exploit the boundedness of $\bar{\eta}^\varepsilon$,
$\bar{v}^\varepsilon$, $\bar{\sigma}^\varepsilon$ and
$\bar{G}^\varepsilon$ to control $\partial_t\bar{\eta}^\varepsilon$,
$\partial_t\bar{v}^\varepsilon$ and to give some convergence
results. The first equation in $(\ref{0503})$ implies that
$\partial_t\bar{\eta}^\varepsilon=\bar{v}^\varepsilon$, therefore
\begin{equation}\label{0515}
\sup_{0\leq
t<T}\|\partial_t\bar{\eta}^\varepsilon(t)\|_{\bar{\mathrm{H}}^3}
=\sup_{0\leq t<T}\|\bar{v}^\varepsilon(t)\|_{\bar{\mathrm{H}}^3}\leq
1.
\end{equation}
By virtue of (\ref{0509}) and (\ref{0514}), the third equation in
(\ref{0503}) yields
\begin{equation}\label{0516}
\lim_{\varepsilon\rightarrow 0}\sup_{0\leq t<T}
\|\mathrm{div}\bar{v}^\varepsilon(t)\|_{\bar{\mathrm{H}}^2}=\lim_{\varepsilon\rightarrow
0}\sup_{0\leq t<T}\|\varepsilon\mathrm{tr}(\bar{G}^\varepsilon
D\bar{v}^\varepsilon )(t)\|_{\bar{\mathrm{H}}^2}= 0.\end{equation}
Expanding the second equation in (\ref{0503}), one sees that
\begin{equation}\begin{aligned}\label{0517}&
\varrho\partial_t\bar{v}^\varepsilon_i
+\partial_i\bar{q}^\varepsilon -\mu\Delta
\bar{v}_i^\varepsilon-B^2\partial_{11}^2\bar{\eta}^\varepsilon_i\\
&={\varepsilon\bar{G}^\varepsilon_{jk}}(\partial_k(\bar{q}^\varepsilon
I_{ij}-\mu((I_{jk}-\varepsilon\bar{G}_{jk}^\varepsilon)\partial_k
\bar{v}^\varepsilon_i+(I_{ik}-\varepsilon\bar{G}_{ik}^\varepsilon)\partial_k
\bar{v}_j^\varepsilon))\\
&\quad -\varepsilon\mu\partial_j(\bar{G}_{jj}^\varepsilon\partial_j
\bar{v}_i^\varepsilon+\bar{G}_{ij}^\varepsilon\partial_j
\bar{v}_j^\varepsilon):=\bar{X}^\varepsilon,\ i=1,2,3,
\end{aligned}\end{equation}
whence, by (\ref{0509}) and (\ref{0514}),
\begin{equation}\label{0518}
\lim_{\varepsilon\rightarrow 0} \sup_{0\leq
t<T}\|\bar{X}^\varepsilon\|_{\bar{\mathrm{H}}^1}=0\end{equation}
and
\begin{equation}\label{0519}
\sup_{0\leq t<T}\|\partial_t\bar{v}^\varepsilon\|_{\bar{\mathrm{H}}^1}\leq
{K_3}\qquad\mbox{for some constant }K_3>0. \end{equation}

By (\ref{0509}), (\ref{0510}), (\ref{0519}) and the sequential
weak-*compactness, we see that up to extraction of a subsequence
(which we still denote using only $\varepsilon$),
\begin{equation}\label{0520}(\bar{\eta}^\varepsilon,\bar{v}^\varepsilon,
\bar{q}^\varepsilon)\rightarrow (\bar{\eta},\bar{v}, \bar{q}) \mbox{
weakly-* in }L^\infty(0,T;\bar{\mathrm{H}}^3(\Omega)),
\end{equation}\begin{equation*}\label{0521}\bar{v}^\varepsilon
\rightarrow \bar{v} \mbox{ weakly-* in }L^\infty(0,T;H^1_0(\Omega))
\end{equation*}
and
\begin{equation}\label{0522}\partial_t\bar{v}^\varepsilon\to\partial_t\bar{v}
 \mbox{ weakly-* in }L^\infty (0,T;\bar{\mathrm{H}}^1(\Omega)).
\end{equation}

 From the lower semi-continuity one gets
\begin{equation}\label{0523}\mathrm{ess}\sup_{0\leq t<T}\| (\bar{\eta},\bar{v},
\bar{q})(t)\|_{\bar{\mathrm{H}}^3 }\leq 1,\;\;\; \mathrm{ess}\sup_{0\leq
t<T}\|\bar{v}\|_{{{H}}^1_0(\Omega)}\leq 1.
\end{equation}
On the other hand, if we use (\ref{0515}), (\ref{0519}), and (\ref{0509}),
 an abstract version of the Arzela-Ascoli theorem (see \cite[Theorem 1.70]{NASII04}), and
the Aubin-Lions Theorem (see \cite[Theorem 1.71]{NASII04}), we obtain that
\begin{eqnarray}\label{0524}&&(\bar{\eta}^\varepsilon,\bar{v}^\varepsilon)
\to (\bar{\eta},\bar{v})\;\mbox{ strongly in } L^{p}(0,T;\bar{\mathrm{H}}^{2}(\Omega))
\;\mbox{ and }C^0([0,T],L^2(\Omega))
\end{eqnarray}
for any $1\leq p<\infty$. This strong convergence, together with the equation
$\partial_t\bar{\eta}^\varepsilon=\bar{v}^\varepsilon$, gives that
\begin{equation}\label{0525}   \begin{array}{l}
                            \partial_t\bar{\eta}^\varepsilon\rightarrow
\partial_t\bar{\eta}\;\;\mbox{ strongly in }L^{p}(0,T;\bar{\mathrm{H}}^{2}(\Omega)).
                     \end{array}
\end{equation}

Rescaling the equations (\ref{0503}), utilizing
(\ref{0516})--(\ref{0518}),  (\ref{0520}), (\ref{0522}) and
(\ref{0525}), the resulting equations thus imply that
\begin{equation}\label{0526} \left\{
                              \begin{array}{ll}
\partial_t \bar{\eta}=\bar{v},\\
\varrho\partial_t\bar{v}+\nabla \bar{\sigma}=\mu\Delta\bar{v}+
B^2\partial_{11}^2\bar{\eta},\\
\mathrm{div}\bar{v}=0
\end{array}
                            \right.
\end{equation}holds a.e. in $(0,T)\times \Omega$.
From (\ref{0524}) and the initial conditions
($\bar{\eta}^\varepsilon,\bar{v}^\varepsilon$)(0)=($\tilde{\eta},\tilde{v}$)(0), it follows that
\begin{equation}\label{0527}
(\bar{\eta}, \bar{v})(0)=(\tilde{\eta},\tilde{v})(0)
\end{equation}
as well.

\subsection{Convergence to the linearized jump conditions}

We proceed to show some convergence results for the jump conditions.
Multiplying (\ref{0503})$_2$ with $\varphi\in
(C_0^\infty((0,T)\times \Omega))^3$, rescaling the resulting
equation, we find that
\begin{equation}\label{0528}\int_0^T\int_{\Omega}(\varrho\partial_t
\bar{v}^\varepsilon-B^2\partial_{11}^2\bar{\eta}^\varepsilon)
\cdot\varphi\mathrm{d}x\mathrm{d}t+\int_0^T\int_{\Omega}(I_{jk}-\varepsilon\bar{G^\varepsilon_{jk}})
\varphi_i\partial_k\bar{T}^\varepsilon_{ij}\mathrm{d}x\mathrm{d}t=0.\end{equation}
Integrating by parts and making use of the second jump condition in
(\ref{0504}), we conclude that
\begin{equation}\begin{aligned}\label{0529}
& \int_0^T\int_{\Omega} \varphi_i\partial_j
\bar{T}_{ij}^\varepsilon\mathrm{d}x\mathrm{d}t +
\int_0^T\int_{\Omega}\bar{T}_{ij}^\varepsilon\partial_j\varphi_i
\mathrm{d}x\mathrm{d}t =-\int_0^T\int_{\mathbb{R}^2}\varphi_i
\llbracket\bar{T}_{i3}^\varepsilon  \rrbracket
\mathrm{d}x'\mathrm{d}t\\
&\quad =-\int_0^T\int_{\mathbb{R}^2}\varphi_i
\llbracket\bar{T}_{ij}^\varepsilon n_j \rrbracket
\mathrm{d}x'\mathrm{d}t+\int_0^T\int_{\mathbb{R}^2}\varphi_i
(\llbracket\bar{T}_{ij}^\varepsilon
n_j\rrbracket-\llbracket\bar{T}_{i3}^\varepsilon \rrbracket)
\mathrm{d}x'\mathrm{d}t  \\
&\quad
=-\int_0^T\int_{\mathbb{R}^2}(g[\varrho]\bar{\eta}_3^\varepsilon
{n}^\varepsilon+{B}^2(e_1+\varepsilon\partial_1\bar{\eta}^\varepsilon)\cdot
n^\varepsilon\llbracket
\partial_1 \bar{\eta}^\varepsilon \rrbracket +\kappa \bar{H}^\varepsilon
n^\varepsilon)\cdot\varphi\mathrm{d}x'\mathrm{d}t\\
&\quad \quad +\int_0^T\int_{\mathbb{R}^2}\varphi_i
(\llbracket\bar{T}_{ij}^\varepsilon
n_j\rrbracket-\llbracket\bar{T}_{i3}^\varepsilon \rrbracket)
\mathrm{d}x'\mathrm{d}t,
\end{aligned}\end{equation}
where $\bar{H}^\varepsilon=H^\varepsilon/\varepsilon$. On the other hand,
\begin{equation}\begin{aligned}\label{0530}
\int_0^T\int_{\Omega}\bar{T}_{ij}^\varepsilon\partial_j\varphi_i
\mathrm{d}x\mathrm{d}t=&\int_0^T\int_\Omega(\bar{q}^\varepsilon\mathrm{div}\varphi-\mu(\nabla
\bar{v}^\varepsilon+\nabla (\bar{v}^{\varepsilon})^T):\nabla
\varphi)\mathrm{d}x\mathrm{d}t\\ &+\varepsilon\int_0^T\int_\Omega\mu
(\bar{G}_{jk}^\varepsilon\partial_k\bar{v}_i^\varepsilon
+\bar{G}_{ik}^\varepsilon\partial_k\bar{v}_j^\varepsilon)\partial_j\varphi_i\mathrm{d}x\mathrm{d}t
\end{aligned}\end{equation}
Putting (\ref{0528})--(\ref{0530}) together, we arrive at
\begin{equation}\begin{aligned}\label{0531}
&\int_0^T\int_{\Omega}(\varrho\partial_t
\bar{v}^\varepsilon-B^2\partial_{11}^2\bar{\eta}^\varepsilon)\cdot\varphi\mathrm{d}x\mathrm{d}t
+ \int_0^T\int_{\Omega}\left(\mu(\nabla \bar{v}^\varepsilon +\nabla
(\bar{v}^{\varepsilon})^T):\nabla
\varphi-\bar{q}^\varepsilon\mathrm{div}\varphi
\right)\mathrm{d}x\mathrm{d}t\\
&\quad \quad-\varepsilon\left[\int_0^T\int_{\Omega}
\varphi_i\bar{G}^\varepsilon_{jk}\partial_k\bar{T}^\varepsilon_{ij}\mathrm{d}x\mathrm{d}t+\int_0^T\int_\Omega\mu
(\bar{G}_{jk}^\varepsilon\partial_k\bar{v}_i^\varepsilon
+\bar{G}_{ik}^\varepsilon\partial_k\bar{v}_j^\varepsilon)\partial_j\varphi_i\mathrm{d}x\mathrm{d}t
\right]\\
&\quad =\int_0^T\int_{\mathbb{R}^2}(g[\varrho]\bar{\eta}_3^\varepsilon
{n}^\varepsilon+{B}^2(e_1+\varepsilon
\partial_1\bar{\eta}^\varepsilon)\cdot n^\varepsilon\llbracket
\partial_1 \bar{\eta}^\varepsilon \rrbracket+\kappa \bar{H}^\varepsilon
n^\varepsilon)\cdot \varphi\mathrm{d}x'\mathrm{d}t\\
&\quad \quad -\int_0^T\int_{\mathbb{R}^2}\varphi_i
(\llbracket\bar{T}_{ij}^\varepsilon
n_j\rrbracket-\llbracket\bar{T}_{i3}^\varepsilon \rrbracket)
\mathrm{d}x'\mathrm{d}t.\end{aligned}\end{equation}

Next, we deal with the limit of the equality (\ref{0531}) as
$\varepsilon\rightarrow 0$. Obviously, by (\ref{0520}),
(\ref{0522}), the bounds (\ref{0509}) and (\ref{0514}),
it is easy to verify that
\begin{equation}\begin{aligned}\label{0532}
&\lim_{\varepsilon\rightarrow
0}\int_0^T\int_{\Omega}(\varrho\partial_t
\bar{v}^\varepsilon-B^2\partial_{11}^2\bar{\eta}^\varepsilon)\cdot
\varphi\mathrm{d}x\mathrm{d}t\\
&\quad +\lim_{\varepsilon\rightarrow 0}
\int_0^T\int_{\Omega}\left(\mu(\nabla \bar{v}^\varepsilon +\nabla
(\bar{v}^{\varepsilon})^T):\nabla
\varphi-\bar{q}^\varepsilon\mathrm{div}\varphi
\right)\mathrm{d}x\mathrm{d}t\\
&\quad =\int_0^T\int_{\Omega}(\varrho\partial_t
\bar{v}-B^2\partial_{11}^2\bar{\eta})\cdot
\varphi\mathrm{d}x\mathrm{d}t+\int_0^T\int_{\Omega}\left(\mu(\nabla
\bar{v} +\nabla \bar{v}^T):\nabla \varphi-\bar{q}\mathrm{div}\varphi
\right)\mathrm{d}x\mathrm{d}t,
\end{aligned}\end{equation}
and
\begin{equation}\label{0533}\lim_{\varepsilon\rightarrow
0}\varepsilon\left[\int_0^T\int_{\Omega}
\varphi_i\bar{G}^\varepsilon_{jk}\partial_k\bar{T}^\varepsilon_{ij}\mathrm{d}x\mathrm{d}t
+\int_0^T\int_\Omega\mu (\bar{G}_{jk}^\varepsilon
\partial_k\bar{v}_i^\varepsilon+\bar{G}_{ik}^\varepsilon\partial_k
\bar{v}_j^\varepsilon)\partial_j\varphi_i\mathrm{d}x\mathrm{d}t
\right]=0.\end{equation}
Thus, it remains to analyze the convergence
of the two integrals on the right hand of  (\ref{0531}).
To this end, using (\ref{0506}) and (\ref{0507}), we rewrite the
first integral as
\begin{equation}\label{0534}\begin{aligned}
&\int_0^T\int_{\mathbb{R}^2}(g[\varrho]\bar{\eta}_3^\varepsilon
{n}^\varepsilon+{B}^2(e_1+\varepsilon
\partial_1\bar{\eta}^\varepsilon)\cdot n^\varepsilon\llbracket
\partial_1 \bar{\eta}^\varepsilon \rrbracket+\kappa \bar{H}^\varepsilon
n^\varepsilon)\cdot \varphi\mathrm{d}x'\mathrm{d}t \\
& =
\int_0^T\int_{\mathbb{R}^2}\bigg\{\frac{g[\varrho]\bar{\eta}_{3}^\varepsilon
e_3}{|e_3+\varepsilon\bar{N}^\varepsilon|}+\frac{\kappa\Delta_{x'}\bar{\eta}_3^\varepsilon
e_3}{(1+\varepsilon \bar{R}^\varepsilon)|e_3+\varepsilon\bar{N}^\varepsilon|^2} \\
&\qquad \qquad\quad \ +\varepsilon\left[\frac{B^2
\bar{F}^\varepsilon
+g[\varrho]\bar{\eta}^\varepsilon_3{\bar{N}^\varepsilon}}{|e_3+\varepsilon\bar{N}^\varepsilon|}
+ \frac{\kappa \bar{L}^\varepsilon}{(1+\varepsilon
\bar{R}^\varepsilon)|e_3+\varepsilon\bar{N}^\varepsilon|^2}\right]\bigg\}\cdot
\varphi\mathrm{d}x'\mathrm{d}t,
\end{aligned}\end{equation}
where
\begin{equation*}\begin{aligned}
\bar{R}^\varepsilon=&2\sum_{i=1}^2\partial_i\bar{\eta}_{i}^\varepsilon+\varepsilon
\sum_{i=1}^2|\partial_i\bar{\eta}^\varepsilon|^2+\varepsilon\left(2\partial_1\bar{\eta}_{1}^\varepsilon+
\varepsilon|\partial_1\bar{\eta}^\varepsilon|^2\right)\left(2\partial_2\bar{\eta}_{2}^\varepsilon+
\varepsilon|\partial_2\bar{\eta}^\varepsilon|^2\right)\\
&-\varepsilon\left(\partial_1\bar{\eta}_2^\varepsilon
+\partial_2\bar{\eta}_1^\varepsilon+\varepsilon\partial_1\bar{\eta}^\varepsilon
\cdot\partial_2\bar{\eta}^\varepsilon\right)^2
,\end{aligned}
\end{equation*}
\begin{equation*}\begin{aligned}
\bar{F}^\varepsilon=&
(\bar{N}_1^\varepsilon+\partial_1\bar{\eta}_{3}^\varepsilon+\varepsilon
\partial_1\bar{\eta}^\varepsilon\cdot \bar{N^\varepsilon})
\llbracket\partial_1\bar{\eta}^\varepsilon\rrbracket,\qquad\qquad\qquad\qquad\qquad\qquad\quad \ \
\end{aligned}\end{equation*}
\begin{equation*}\begin{aligned}&\bar{L}^\varepsilon=\{(2\partial_1\bar{\eta}_{1}^\varepsilon+
\varepsilon|\partial_1\bar{\eta}^\varepsilon|^2)\partial_2^2\bar{\eta}^\varepsilon_3+
(2\partial_2\bar{\eta}_{2}^\varepsilon+
\varepsilon|\partial_2\bar{\eta}^\varepsilon|^2)\partial_1^2\bar{\eta}_3^\varepsilon
\\&\qquad \ -
2(\partial_1\bar{\eta}_{2}^\varepsilon+\partial_2\bar{\eta}_{1}^\varepsilon
+\varepsilon\partial_1\bar{\eta}^\varepsilon\cdot
\partial_2\bar{\eta}^\varepsilon)\partial_1\partial_2\bar{\eta}^\varepsilon_3\\
&\qquad \
+\bar{\mathcal{H}}^\varepsilon\cdot\bar{N^\varepsilon}](e_3+\bar{N^\varepsilon})
+\Delta_{x'}\bar{\eta}_3^\varepsilon\bar{N}^\varepsilon,
\quad\qquad\qquad\qquad\qquad\qquad \qquad
\end{aligned}
\end{equation*}and
\begin{equation*}\mathcal{\bar{H}}^\varepsilon={|e_1+\varepsilon\partial_1
\bar{\eta}^\varepsilon|^2\partial_2^2\bar{\eta}^\varepsilon
-2(e_1+\varepsilon\partial_1\bar{\eta}^\varepsilon)\cdot(e_2+\varepsilon\partial_2
\bar{\eta}^\varepsilon)\partial_1\partial_2\bar{\eta}^\varepsilon+
|e_2+\varepsilon\partial_2\bar{\eta}^\varepsilon|^2\partial_1^2
\bar{\eta}^\varepsilon}.
\end{equation*}

 Clearly, by the trace theorem, (\ref{0509}) and (\ref{050911}) in terms of
 $\bar{\eta}^\varepsilon$, $\sup_{0\leq t<T}\|(\bar{R}^\varepsilon,
\bar{F}^\varepsilon,\bar{N}^\varepsilon)(t)\|_{L^\infty(\mathbb{R}^2)}$
and $\sup_{0\leq t<T}\|\bar{L}^\varepsilon(t)\|_{L^2(\mathbb{R}^2)}$
are uniformly bounded. Hence, as $\varepsilon\rightarrow 0$, we have
$1+\varepsilon\bar{R}^\varepsilon>0$,
$|e_3+\varepsilon\bar{N}^\varepsilon|>0$, and
\begin{equation}\label{0535}\varepsilon\left[\frac{B^2
\bar{F}^\varepsilon
+g[\varrho]\bar{\eta}^\varepsilon_3{\bar{N}^\varepsilon}}{|e_3+\varepsilon\bar{N}^\varepsilon|}
+ \frac{\kappa \bar{L}^\varepsilon}{(1+\varepsilon
\bar{R}^\varepsilon)|e_3+\varepsilon\bar{N}^\varepsilon|^2}\right]\rightarrow
0\mbox{ strongly in }L^\infty (0,T; L^2(\mathbb{R}^2)).
\end{equation}
Moreover, by virtue of (\ref{0524}), (\ref{0509}) and  (\ref{050911})
in terms of $\bar{\eta}^\varepsilon$, using integrating by parts, the trace theorem
and dominated convergence theorem, we deduce that
\begin{equation}\label{0536}\begin{aligned}&\lim_{\varepsilon\rightarrow 0}
\int_0^T\int_{\mathbb{R}^2}\frac{\Delta_{x'}\bar{\eta}_3^\varepsilon
e_3\cdot \varphi}{(1+\varepsilon
\bar{R}^\varepsilon)|e_3+\varepsilon\bar{N}^\varepsilon|^2}\mathrm{d}x'\mathrm{d}t\\
&=-\sum_{i=1}^2\lim_{\varepsilon\rightarrow
0}\int_0^T\int_{\mathbb{R}^2}\frac{\partial_i\bar{\eta}^\varepsilon_3\partial_i
\varphi_3}{(1+\varepsilon
\bar{R}^\varepsilon)|e_3+\varepsilon\bar{N}^\varepsilon|^2}\mathrm{d}x'\mathrm{d}t\\
&\quad -\sum_{i=1}^2\lim_{\varepsilon\rightarrow
0}\int_0^T\int_{\mathbb{R}^2}\varphi_i\partial
_i\bar{\eta}^\varepsilon_3\partial_i\left[\frac{1}{(1+\varepsilon
\bar{R}^\varepsilon)|e_3+\varepsilon\bar{N}^\varepsilon|^2}\right]
 \mathrm{d}x'\mathrm{d}t   \\
&=- \sum_{i=1}^2\int_0^T\int_{\mathbb{R}^2}\partial_i\bar{\eta}_3\partial_i
\varphi_3 \mathrm{d}x'\mathrm{d}t
\end{aligned}\end{equation}
and
\begin{equation}\label{0537}\begin{aligned}
\int_0^T\int_{\mathbb{R}^2}\frac{g[\varrho]\bar{\eta}_{3}^\varepsilon
e_3 \cdot
\varphi}{|e_3+\varepsilon\bar{N}^\varepsilon|}\mathrm{d}x'\mathrm{d}t
\rightarrow
\int_0^T\int_{\mathbb{R}^2}g[\varrho]\bar{\eta}_{3}\varphi_3\mathrm{d}x'\mathrm{d}t.
\end{aligned}\end{equation}
In view of (\ref{0535})--(\ref{0537}) and (\ref{0534}), we find that
\begin{equation}\begin{aligned}\label{0538}
&\lim_{\varepsilon\rightarrow
0}\int_0^T\int_{\mathbb{R}^2}(g[\varrho]\bar{\eta}_3^\varepsilon
{n}^\varepsilon+{B}^2(e_1+\varepsilon
\partial_1\bar{\eta}^\varepsilon)\cdot n^\varepsilon\llbracket
\partial_1 \bar{\eta}_i^\varepsilon \rrbracket+\kappa \bar{H}^\varepsilon
n^\varepsilon)\cdot \varphi\mathrm{d}x'\mathrm{d}t  \\
& \quad =\int_0^T\int_{\mathbb{R}^2}(g[\varrho]\bar{\eta}_{3}
\varphi_3-\kappa\sum_{i=1}^2\partial_{i}\bar{\eta}_3)
\partial_i\varphi_3\mathrm{d}x'\mathrm{d}t.
\end{aligned}\end{equation}

Finally, using (\ref{050911}), we easily verify that
\begin{equation}\label{0539}
\lim_{\varepsilon\rightarrow 0}\int_0^T\int_{\mathbb{R}^2}\varphi_i
(\llbracket\bar{T}_{ij}^\varepsilon
n_j\rrbracket-\llbracket\bar{T}_{i3}^\varepsilon \rrbracket)
\mathrm{d}x'\mathrm{d}t=0.\end{equation}
Consequently, thanks to (\ref{0532}), (\ref{0533}), (\ref{0538}) and (\ref{0539}),
we take to the limit in (\ref{0531}) to arrive at
\begin{equation}\begin{aligned}\label{0540}
&\int_0^T\int_{\Omega}(\varrho\partial_t
\bar{v}-B^2\partial_{11}^2\bar{\eta})\varphi\mathrm{d}x\mathrm{d}t+
\int_0^T\int_{\Omega}\left(\mu(\nabla \bar{v} +\nabla
\bar{v}^T):\nabla \varphi-\bar{q}\mathrm{div}\varphi
\right)\mathrm{d}x\mathrm{d}t \\
&\quad =\int_0^T\int_{\mathbb{R}^2}(g[\varrho]\bar{\eta}_{3}e_3-\kappa\sum_{i=1}^2
\partial_{i}\bar{\eta}_3 \partial_i\varphi_3)\mathrm{d}x'\mathrm{d}t.
\end{aligned}\end{equation}

On the other hand, similarly to (\ref{0423}), we multiply
(\ref{0526})$_2$ by $\varphi$ and integrate over $(0,T)\times\Omega$ to obtain
\begin{equation}\begin{aligned}\label{0541}
&\int_0^T\int_{\Omega}(\varrho\partial_t
\bar{v}-B^2\partial_{11}^2\bar{\eta})\cdot\varphi\mathrm{d}x\mathrm{d}t+
\int_0^T\int_{\Omega}\left(\mu(\nabla \bar{v} +\nabla
\bar{v}^T):\nabla \varphi-\bar{q}\mathrm{div}\varphi
\right)\mathrm{d}x\mathrm{d}t\\
&=\int_0^T\int_{\mathbb{R}^2}\Big( \bar{q}_+I-\mu_+(\nabla
\bar{v}_++\nabla \bar{v}^T_+)-\big(\bar{q}_-I-\mu_-(\nabla
\bar{v}_-+\nabla \bar{v}_-^T)\big) \Big)
e_3\cdot\varphi\mathrm{d}x'\mathrm{d}t ,
\end{aligned}\end{equation}
where we have used the fact that $\gamma_\pm (q)=q_\pm$ at the interface
$\{x_3=0\}$ since $q\in \mathrm{\bar{H}}^3(\Omega)$ for a.e. $t\in (0,T)$.
Comparing (\ref{0540}) with (\ref{0541}), we conclude that
\begin{eqnarray}
&& \int_{\mathbb{R}^2}\Big( \bar{q}_+I-\mu_+(\nabla \bar{v}_++\nabla
\bar{v}^T_+)-\big( \bar{q}_-I-\mu_-(\nabla \bar{v}_-+\nabla
\bar{v}_-^T)\big)\Big) e_3\cdot \phi\mathrm{d}x' \nonumber  \\
&& \quad =\int_{\mathbb{R}^2} \Big( g[\varrho]\bar{\eta}_3
\phi_3-\kappa\sum_{i=1}^2\partial_i\bar{\eta}_{3}\partial_i\phi_3\Big)\mathrm{d}x'
\;\;\mbox{ for a.e. }t\in (0,T)\label{0542}\end{eqnarray}
holds for any $\phi\in H_0^1(\mathbb{R}^2)$.

\subsection{Contradiction argument}
Notice that $\bar{\eta}_3$ in (\ref{0542})  includes both cases of
$\bar{\eta}_{+3}$  and $\bar{\eta}_{-3}$. In view of Definition
\ref{def:0401}, we find that $(\bar{\eta},\bar{v}, \bar{q})$ is just
a strong solution of the linearized problem
(\ref{0134})--(\ref{0136}). By Remark \ref{rem:0402},
$(\tilde{\eta},\tilde{v},\tilde{q})$ constructed in Theorem
\ref{thm:0201} is also a strong solution of
(\ref{0134})--(\ref{0136}). Moreover,
$\tilde{\eta}(0)=\bar{\eta}(0)$ and $\tilde{v}(0)=\bar{v}(0)$ (see
(\ref{0527})). Then, according to Theorem \ref{thm:0401},
\begin{equation*}
\bar{v}=\tilde{v}\mbox{ on }[0,{T})\times \Omega.
\end{equation*}
Hence, we may chain together the inequalities (\ref{0523}) and
(\ref{0502}) to get
\begin{equation*}
2\leq \sup_{{T}/2\leq t<{T}}\|
\tilde{v}(t)\|_{\bar{\mathrm{H}}^3}\leq \sup_{0\leq t<{T}}\|
\bar{v}\|_{\bar{\mathrm{H}}^3}\leq 1,
\end{equation*}
which is a contraction. Therefore, the perturbed problem does not
have the global stability of order $k$ for any $k\geq 3$. This
completes the proof of Theorem \ref{thm:0202}.
\\[3mm]
 {\bf{Acknowledgements.}}\quad
 The research of Fei Jiang was supported by the Fujian
Provincial Department of Science and Technology (Grant No.JK2009045)
and NSFC (Grant No.11101044), and the research of Song Jiang by the
National Basic Research Program under the Grant 2011CB309705 and
NSFC (Grant No.40890154).


\bibliographystyle{model1-num-names}







\end{document}